\documentstyle[11pt]{article}

\font\tenfrak=eufm10
\font\sevenfrak=eufm7
\font\fivefrak=eufm5
\newfam\frakfam \def\frak{\fam\frakfam\tenfrak} \textfont\frakfam=\tenfrak
  \scriptfont\frakfam=\sevenfrak  \scriptscriptfont\frakfam=\fivefrak

\begin{document}
    \pagestyle{plain}
    \setlength{\baselineskip}{1.2\baselineskip}
    \setlength{\parindent}{\parindent}

\title{{\bf Poisson harmonic forms, Kostant harmonic forms, and 
the $S^1$-equivariant cohomology of $K/T$}}
\author{Sam Evens\thanks{Research partially supported by NSF 
grant DMS-9623322;}
\hspace{.04in}
 and Jiang-Hua Lu
\thanks{Research partially supported by
an NSF Postdoctorial Fellowship.} \\
Department of Mathematics, University of Arizona, Tucson, AZ 85721\\
evens@math.arizona.edu, jhlu@math.arizona.edu}
\maketitle

\newtheorem{thm}{Theorem}[section]
\newtheorem{lem}[thm]{Lemma}
\newtheorem{prop}[thm]{Proposition}
\newtheorem{cor}[thm]{Corollary}
\newtheorem{rem}[thm]{Remark}
\newtheorem{exam}[thm]{Example}
\newtheorem{nota}[thm]{Notation}
\newtheorem{dfn}[thm]{Definition}
\newtheorem{ques}[thm]{Question}
\newtheorem{eq}{thm}

\newcommand{\rw}{\rightarrow}
\newcommand{\lrw}{\longrightarrow}
\newcommand{\rhu}{\rightharpoonup}
\newcommand{\lhu}{\leftharpoonup}
\newcommand{\Map}{\longmapsto}
\newcommand{\qed}{\begin{flushright} {\bf Q.E.D.}\ \ \ \ \
                  \end{flushright} }
\newcommand{\beqa}{\begin{eqnarray*}}
\newcommand{\eeqa}{\end{eqnarray*}}

\newcommand{\la}{\mbox{$\langle$}}
\newcommand{\ra}{\mbox{$\rangle$}}
\newcommand{\ot}{\mbox{$\otimes$}}
\newcommand{\xa}{\mbox{$x_{(1)}$}}
\newcommand{\xb}{\mbox{$x_{(2)}$}}
\newcommand{\xc}{\mbox{$x_{(3)}$}}
\newcommand{\ya}{\mbox{$y_{(1)}$}}
\newcommand{\yb}{\mbox{$y_{(2)}$}}
\newcommand{\yc}{\mbox{$y_{(3)}$}}
\newcommand{\yd}{\mbox{$y_{(4)}$}}
\renewcommand{\aa}{\mbox{$a_{(1)}$}}
\newcommand{\ab}{\mbox{$a_{(2)}$}}
\newcommand{\ac}{\mbox{$a_{(3)}$}}
\newcommand{\ad}{\mbox{$a_{(4)}$}}
\newcommand{\ba}{\mbox{$b_{(1)}$}}
\newcommand{\bt}{\mbox{$b_{(2)}$}}
\newcommand{\bc}{\mbox{$b_{(3)}$}}
\newcommand{\ca}{\mbox{$c_{(1)}$}}
\newcommand{\cb}{\mbox{$c_{(2)}$}}
\newcommand{\cc}{\mbox{$c_{(3)}$}}
\newcommand{\uo}{\mbox{$\underbar{o}$}}
\newcommand{\eo}{\mbox{$e$}}

\newcommand{\ts}{\mbox{$\sigma$}}
\newcommand{\las}{\mbox{${}_{\sigma}\!A$}}
\newcommand{\lasone}{\mbox{${}_{\sigma'}\!A$}}
\newcommand{\ras}{\mbox{$A_{\sigma}$}}
\newcommand{\rds}{\mbox{$\cdot_{\sigma}$}}
\newcommand{\lds}{\mbox{${}_{\sigma}\!\cdot$}}

\newcommand{\bb}{\mbox{$\bar{\beta}$}}
\newcommand{\bg}{\mbox{$\bar{\gamma}$}}

\newcommand{\id}{\mbox{${\rm id}$}}
\newcommand{\Fun}{\mbox{${\rm Fun}$}}
\newcommand{\End}{\mbox{${\rm End}$}}
\newcommand{\Hom}{\mbox{${\rm Hom}$}}
\newcommand{\Ker}{\mbox{${\rm Ker}$}}
\renewcommand{\Im}{\mbox{${\rm Im}$}}

\newcommand{\ta}{\mbox{${\mbox{$\scriptscriptstyle A$}}$}}
\newcommand{\ms}{\mbox{${\mbox{$\scriptscriptstyle M$}}$}}
\newcommand{\ap}{\mbox{$A_{\mbox{$\scriptscriptstyle P$}}$}}
\newcommand{\tx}{\mbox{$\mbox{$\scriptscriptstyle X$}$}}
\newcommand{\tE}{\mbox{$\mbox{$\scriptscriptstyle E$}$}}
\newcommand{\ty}{\mbox{$\mbox{$\scriptscriptstyle Y$}$}}
\newcommand{\kt}{\mbox{$K_{\tx}$}}
\newcommand{\pk}{\mbox{$\pi_{\lambda}^{\scriptscriptstyle K_{X}}$}}
\newcommand{\kk}{\mbox{$K \times_{\scriptscriptstyle K_{X}} (K_{\tx}/T)$}}

\newcommand{\wm}{W^{\tx}}
\newcommand{\wx}{W_{\tx}}
\newcommand{\dw}{\dot{w}}
\newcommand{\fw}{f^{w}_{w_1}}
\newcommand{\lw}{C^{w}_{w_1}}

\newcommand{\sww}{s^{w}_{w_1}}
\newcommand{\sw}{\Sigma_w}
\newcommand{\swa}{\Sigma_{w_1}}
\newcommand{\swb}{\Sigma_{w_1 w_2}}
\newcommand{\cw}{C_{\dot{w}}}
\newcommand{\cwa}{C_{\dot{w_1}}}

\newcommand{\pp}{\mbox{$\pi_{\mbox{$\scriptscriptstyle P$}}$}}
\newcommand{\pg}{\mbox{$\pi_{\mbox{$\scriptscriptstyle G$}}$}}
\newcommand{\pge}{\mbox{$\pi_{\Delta_{+}}$}}
\newcommand{\iep}{\mbox{$i_{\exp_{\wedge} \pi}$}}
\newcommand{\ienp}{\mbox{$i_{\exp_{\wedge}(- \pi)}$}}

\newcommand{\Xa}{\mbox{$X_{\alpha}$}}
\newcommand{\Ya}{\mbox{$Y_{\alpha}$}}
\newcommand{\pix}{\mbox{$\pi_{\tx,\lambda}$}}
\newcommand{\fix}{\mbox{$\fl_{\tx,\lambda}$}}
\newcommand{\fll}{\mbox{$\fl_{\lambda}$}}

\newcommand{\bix}{\mbox{$b_{\tx,\lambda}$}}
\newcommand{\dix}{\mbox{$d_{\tx,\lambda}$}}
\newcommand{\pxl}{\mbox{$\partial_{\tx, \lambda}$}}
\newcommand{\tml}{\mbox{$\tilde{m}_{\lambda}$}}
\newcommand{\tse}{\mbox{$T^{*}_{e}(K/T)$}}
\newcommand{\te}{\mbox{$T_{e}(K/T)$}}
\newcommand{\pinf}{\mbox{$\pi_{\infty}$}}
\newcommand{\ponf}{\mbox{$\partial_{\infty}$}}
\newcommand{\ii}{\mbox{$I_{\infty}$}}
\newcommand{\iis}{\mbox{$I_{\infty}^{*}$}}
\newcommand{\il}{\mbox{$I_{\lambda}$}}
\newcommand{\ils}{\mbox{$I_{\lambda}^{*}$}}
\newcommand{\pl}{\mbox{$\pi_{\lambda}$}}
\newcommand{\pal}{\mbox{$\partial_{\lambda}$}}

\newcommand{\hus}{\mbox{$H_{S^1}(K/T, \C)$}}
\newcommand{\hls}{\mbox{$H^{S^1}(K/T, \C)$}}
\newcommand{\cda}{\mbox{$\sigma^{(w_1)}$}}
\newcommand{\cdb}{\mbox{$\sigma^{(w_2)}$}}
\newcommand{\cd}{\mbox{$\sigma^{(w)}$}}

\newcommand{\dfl}{\mbox{$d_{{\frak l}_{\lambda}} $}}
\newcommand{\bfl}{\mbox{$b_{{\frak l}_{\lambda}} $}}

\newcommand{\el}{\mbox{$e^{- \lambda}$}}
\newcommand{\pkl}{\mbox{$\pi^{{\scriptscriptstyle K_{X}}}_{\lambda}$}}

\newcommand{\dpi}{\mbox{$\delta_{\pi}$}}

\newcommand{\asemi}{\mbox{$\ap \#_{\sigma} A^*$}}
\newcommand{\dsemi}{\mbox{$A \#_{\Delta} A^*$}}

\newcommand{\semi}{\mbox{$\times_{{\frac{1}{2}}}$}}
\newcommand{\fd}{\mbox{${\frak d}$}}
\newcommand{\fa}{\mbox{${\frak a}$}}
\newcommand{\ft}{\mbox{${\frak t}$}}
\newcommand{\fk}{\mbox{${\frak k}$}}
\newcommand{\fg}{\mbox{${\frak g}$}}
\newcommand{\fq}{\mbox{${\frak q}$}}
\newcommand{\fl}{\mbox{${\frak l}$}}
\newcommand{\fs}{\mbox{${\frak s}$}}

\newcommand{\flp}{\mbox{${\frak l}_{p}$}}
\newcommand{\fh}{\mbox{${\frak h}$}}
\newcommand{\fn}{\mbox{${\frak n}$}}
\newcommand{\bfny}{\mbox{${\bar{\frak n}}_{\ty}$}}

\newcommand{\fp}{\mbox{${\frak p}$}}
\newcommand{\fb}{\mbox{${\frak b}$}}
\newcommand{\fbp}{\mbox{${\frak b}_{+}$}}
\newcommand{\fbm}{\mbox{${\frak b}_{-}$}}
\newcommand{\fnp}{\mbox{${\frak n}_{+}$}}
\newcommand{\fnm}{\mbox{${\frak n}_{-}$}}
\newcommand{\fgs}{\mbox{${\frak g}^*$}}
\newcommand{\wg}{\mbox{$\wedge {\frak g}$}}
\newcommand{\wgs}{\mbox{$\wedge {\frak g}^*$}}
\newcommand{\wxl}{\mbox{$x_1 \wedge x_2 \wedge \cdots \wedge x_l$}}
\newcommand{\wxk}{\mbox{$x_1 \wedge x_2 \wedge \cdots \wedge x_k$}}
\newcommand{\wyl}{\mbox{$y_1 \wedge y_2 \wedge \cdots \wedge y_l$}}
\newcommand{\wxkm}{\mbox{$x_1 \wedge x_2 \wedge \cdots \wedge x_{k-1}$}}
\newcommand{\wxik}{\mbox{$\xi_1 \wedge \xi_2 \wedge \cdots \wedge \xi_k$}}
\newcommand{\wxikm}{\mbox{$\xi_1 \wedge \cdots \wedge \xi_{k-1}$}}
\newcommand{\wetal}{\mbox{$\eta_1 \wedge \eta_2 \wedge \cdots \wedge \eta_l$}}

\newcommand{\winv}{\mbox{$(\wedge \fg_{1}^{\perp})^{\fg_1}$}}
\newcommand{\wetak}{\mbox{$\eta_1 \wedge \cdots \wedge \eta_k$}}
\newcommand{\gonep}{\mbox{$\fg_{1}^{\perp}$}}
\newcommand{\wonep}{\mbox{$\wedge \fg_{1}^{\perp}$}}

\newcommand{\cala}{\mbox{${\cal A}$}}
\newcommand{\calv}{\mbox{${\cal V}$}}
\newcommand{\pdp}{\mbox{$\partial_{\pi}$}}

\newcommand{\db}{\mbox{$\fd = \fg \bowtie \fgs$}}
\newcommand{\fds}{\mbox{${\scriptscriptstyle {\frak d}}$}}

\newcommand{\Gs}{\mbox{$G^*$}}
\newcommand{\pis}{\mbox{$\pi_{\sigma}$}}
\newcommand{\ea}{\mbox{$E_{\alpha}$}}
\newcommand{\eb}{\mbox{$E_{-\alpha}$}}
\newcommand{\Bm}{\mbox{$ {}^B \! M$}}
\newcommand{\kBm}{\mbox{$ {}^B \! M^k$}}
\newcommand{\Bb}{\mbox{$ {}^B \! b$}}
\newcommand{\epe}{\mbox{$\epsilon$}}
\newcommand{\eot}{\mbox{${\epe \over 2}$}}

\newcommand{\cfg}{\mbox{$C(\fg \oplus \fgs)$}}
\newcommand{\ps}{\mbox{$\pi^{\#}$}}
\newcommand{\backl}{\mathbin{\vrule width1.5ex height.4pt\vrule height1.5ex}}
 
\newcommand{\bx}{\mbox{${\bar{x}}$}}
\newcommand{\by}{\mbox{${\bar{y}}$}}
\newcommand{\bz}{\mbox{${\bar{z}}$}}
\newcommand{\pgs}{\mbox{${\pi_{\mbox{\tiny G}^{*}}}$}}

\newcommand{\tlp}{\mbox{$\tilde{\pi}$}}
\newcommand{\tp}{\mbox{$\tilde{\pi}$}}
\newcommand{\sn}{\mbox{$s_{\scriptscriptstyle N}$}}
\newcommand{\tn}{\mbox{$t_{\scriptscriptstyle N}$}}
\newcommand{\sm}{\mbox{$s_{\scriptscriptstyle M}$}}
\newcommand{\tm}{\mbox{$t_{\scriptscriptstyle M}$}}
\newcommand{\en}{\mbox{$\epsilon_{\scriptscriptstyle N}$}}
\newcommand{\mem}{\mbox{$\epsilon_{\scriptscriptstyle M}$}}

\newcommand{\C}{\mbox{${\Bbb C}$}}
\newcommand{\Z}{\mbox{${\Bbb Z}$}}

\newcommand{\R}{\mbox{${\Bbb R}$}}
\renewcommand{\a}{\mbox{$\alpha$}}

\input amssym.def
\input amssym.tex

\begin{abstract}
We characterize the harmonic forms on a flag manifold $K/T$
defined by Kostant in 1963 in terms of a Poisson structure.
Namely, they are ``Poisson harmonic" with respect
to the so-called Bruhat Poisson structure on $K/T$. 
This enables us to give Poisson geometrical proofs of many of the
special properties of these harmonic forms. In particular,
we construct explicit
representatives for the Schubert basis of the $S^1$-equivariant cohomology
of $K/T$, where the $S^1$-action is defined by $\rho$. 
Using a simple argument in equivariant cohomology, 
we recover 
the connection between the Kostant harmonic forms and 
the Schubert calculus on $K/T$ that was found by Kostant and Kumar in 1986.
We also show that the Kostant harmonic forms are limits of the more
familiar Hodge harmonic forms with respect to a family of Hermitian
metrics.
\end{abstract}

\tableofcontents

\section{Introduction}
\label{sec_intro}

Let $K$ be a compact semi-simple Lie group and $T \subset K$
a maximal torus.
In \cite{ko:63}, Kostant introduced a degree $-1$ operator 
$\partial$ on the
space of $K$-invariant complex valued differential
forms on the flag manifold $K/T$. He then constructed certain 
forms $s^w$ on $K/T$, for $w \in W$, the Weyl group of 
$(K,T)$, that are $(d, \partial)$-harmonic, where
$d$ is the de Rham differential. 
In this work, we establish 
connections between these harmonic forms and 
the geometry of certain Poisson structures on $K/T$.

\bigskip
We think it is fair to say that, although 
introduced more than 30 years ago, the Kostant harmonic forms
have remained mysterious to this day.
What is especially mysterious is the fact that they possess many
special properties (to be reviewed later in this introduction)
that other kinds of harmonic forms do not. As was emphasized
in \cite{ko:63},
this is all due to the fact that the forms $s^w$ are $\partial$-closed.
Thus a  conceptual understanding of the operator $\partial$
is the key to unveil the mystery about the Kostant harmonic forms.
We point out that $\partial$ is in general not
the adjoint operator of $d$ with respect to any Hermitian
metric on $K/T$, so the harmonicity of the forms $s^w$
is not the same as that in the sense of Hodge.

\bigskip
In this paper, we characterize 
the Kostant operator $\partial$, and thus the Kostant
harmonic forms as well, using the
so-called Bruhat Poisson structure on $K/T$ \cite{lu-we:poi}.
We also show that $\partial$ is the limit of the adjoint 
operators of $d$ with respect to a family of 
Hermitian metrics that come
from a family of symplectic structures on $K/T$.

\bigskip
Given an orientable manifold
$P$ with a Poisson structure $\pi$ and a volume form $\mu$, we
consider the degree $-1$ {\bf Koszul-Brylinski operator} 
$\partial_{\pi, \mu}$
on the space of differential forms on $P$ defined by
\[
\partial_{\pi, \mu} \,  =  \, i_{\pi} d \, - \, d i_{\pi}
\, + \, i_{\theta_{\mu}},
\]
where $\theta_{\mu}$ is the modular vector field of $\pi$ with
respect to $\mu$, and $i_{\pi}$ and $i_{\theta_{\mu}}$
are the contraction operators defined by $\pi$ and $\theta_{\mu}$
respectively (see Section \ref{sec_poi-harm} for details). 
We say that
a form $\xi$ on $P$ is ``{\bf Poisson harmonic with respect
to $\pi$ and $\mu$}" if $d \xi = \partial_{\pi, \mu} \xi = 0$.
This is a modified version of a notion
introduced by Brylinski in \cite{by:poi}.

\bigskip
The Bruhat Poisson structure on $K/T$, first introduced in
\cite{lu-we:poi} (see also \cite{soi:compact}),
 has its origin from quantum groups. 
It is so named because its symplectic
leaves are exactly all the Bruhat (or Schubert) cells in
$K/T$.
In this paper, the Bruhat Poisson structure
is denoted by $\pinf$ (for reasons given below)
and its Koszul-Brylinski operator (using a $K$-invariant volume form) 
by $\ponf$.  Our Theorem \ref{thm_main-1}
says that the Kostant operator $\partial$
is related to $\ponf$ by
$\partial = J \ponf J^{-1}$ on the
space $C$  of $K$-invariant complex valued 
differential forms on $K/T$, where $J$ 
is a standard complex structure on $K/T$. Consequently
(Corollary \ref{cor_main-2}), the Kostant harmonic forms
$s^w$, for $w \in W$, which have pure bi-degree with 
respect to the bi-grading on $C$ defined by $J$,
are Poisson harmonic with respect to the Bruhat Poisson structure.
In fact, Theorem 4.5 in \cite{ko:63} can be reformulated as saying
that every de Rham cohomology class of $K/T$ has 
a unique $K$-invariant representative that is Poisson harmonic.
Once this is proved, we show that
the special properties of these Kostant harmonic forms follow
from fairly general arguments in Poisson geometry.
In particular, we construct 
from the forms $s^w, w \in W,$ explicit 
representatives for the Schubert basis of the
$S^1$-equivariant cohomology of $K/T$, where the 
$S^1$ action is defined by $\rho$, half of the sum of all
positive roots (for a choice of such roots). 
Using a simple argument in equivariant cohomology, we then
show geometrically how the Kostant harmonic forms can be used to
describe the ring structure on the de Rham cohomology
$H^{\bullet}(K/T, \C)$ of $K/T$, a fact first proved by Kostant and 
Kumar in \cite{k-k:integral}. To make connections between the 
Kostant harmonic forms and Hodge harmonic forms, 
we employ a family of symplectic structures $\pi_{\lambda}$
on $K/T$ for $\lambda \in \ft^*$ regular, where $\ft$ 
is the Lie algebra of $T$, which has the
property that $\pi_{\lambda}$ goes to the
Bruhat Poisson structure when $\lambda \rightarrow \infty$
within the positive Weyl chamber (and thus the notation
$\pinf$ for the latter). Each $\pi_{\lambda}$ gives rise to 
a Hermitian metric $h_{\lambda}$ on $K/T$. Using the
Koszul-Brylinski operators for the  $\pi_{\lambda}$'s, we show that
the Kostant operator $\partial$ is the limit of
the adjoint operators of $d$ with respect to the Hermitian
metrics $h_{\lambda}$ as $\lambda \rightarrow \infty$.
Correspondingly, the Kostant harmonic forms are shown to 
be limits of ordinary Hodge harmonic forms.

\bigskip
To further explain the content of this paper, we
recall Kostant's construction of the harmonic forms
in \cite{ko:63}. In fact, we
will recall the main results in \cite{ko:63}.

\bigskip
Since $K$ is compact, the (complex-valued) de Rham cohomology
$H^{\bullet}(K/T, \C)$ of $K/T$ can be calculated from the space
$C = \Omega(K/T, \C)^K$ of $K$-invariant
complex valued differential forms on $K/T$. Denote by $\fk$ and $\ft$
the Lie algebras of $K$ and $T$ respectively and by $\fg = \fk_{\Bbb C}$
the complexification of $\fk$.
Then
we can use the Killing form of $\fg$ to identify the 
complexified cotangent space of $K/T$ at the base point $e$
with the vector space $\fn_{-} \oplus \fn_{+}$, where $\fn_{-}$ and
$\fn_{+}$ are the Lie subalgebras of $\fg $
 spanned respectively by the positive and negative root vectors. 
Consequently, we have the identification (see (\ref{eq_I}))
\begin{equation}
\label{eq_i-intro}
I: \, \, (\wedge^{\bullet} (\fn_{-} \oplus \fn_{+}))^T \, 
\stackrel{\sim}{\lrw}
C^{\bullet} \, = \, \Omega^{\bullet}(K/T)^K.
\end{equation}
Now equip $\fn_{-} \oplus \fn_{+}$ with the direct sum Lie 
algebra structure, where $\fn_{-}$ and $\fn_{+}$ have the 
Lie subalgebra structures of $\fg$.
Let $b_{{\frak n}_{-} \oplus {\frak n}_{+}}$ be
the Chevalley-Eilenberg boundary operator for this Lie algebra.
Then Kostant introduced  in \cite{ko:63}
the degree $-1$ operator $\partial$ on
$C = \Omega(K/T, \C)^K$ by
\[
\partial \, := \,  - \,I b_{{\frak n}_{-} \oplus {\frak n}_{+}} I^{-1}: \, \,
C^q \lrw C^{q-1}.
\]
The first main theorem in \cite{ko:63} says that the two operators
$d$ and $\partial$, where $d$ is the de Rham 
differential, are ``disjoint"
in the sense that $\Im (d) \cap \Ker (\partial)
= \Im (\partial) \cap \Ker (d) = 0$. Set $S = d \partial + 
\partial d$ and call it the ``Laplacian" of $d$ and $\partial$.
Then it follows immediately from the disjointness of $d$ and $\partial$
that $\Ker (S) = \Ker (d) \cap \Ker (\partial)$ and that
the natural maps
\beqa
\psi_{d, S}: & & \Ker (S) \lrw H(C, d): \, \, \xi \Map [\xi]_d \\
\psi_{\partial, S}: & &  \Ker (S) \lrw H(C, \partial): \, \, 
\xi \Map [\xi]_\partial
\eeqa
are isomorphisms of graded vector spaces. Elements in $\Ker(S)$
are said to be $(d, \partial)$-harmonic or simply harmonic
in \cite{ko:63} \footnote{To 
distinguish from the notion of Poisson harmonic forms,
we will say throughout this paper
  that such forms are 
``Kostant harmonic" or ``harmonic in the sense of Kostant"}.
Our Example \ref{exam_not-adjoint} shows that
the operator $S$ is in general not semi-simple, so
$\partial$ is in general not the 
adjoint of $d$ with respect to any Hermitian metric
on $K/T$.

\bigskip
What is done next in
\cite{ko:63} is the construction of
 a special basis of $\Ker (S)$ which we explain now.
The Lie algebra homology $H_{\bullet}(\fn_{+})$
of $\fn_{+}$ as a $T$-module had been determined 
earlier by Kostant in \cite{ko:61} to have weights exactly of the
form $\rho - w \rho$, for $w \in W$, and with multiplicity one
for each weight. By viewing
$\fn_{-}$ as the contragradient of $\fn_{+}$ as $T$-modules via
the Killing form of $\fg$, we see by Schur's Lemma that there is a 
canonical basis for the Lie algebra homology
$(H_{\bullet}(\fn_{-} \oplus \fn_{+}))^T \cong
(H_{\bullet}(\fn_{-}) \otimes H_{\bullet}
(\fn_{+}))^T$ 
and thus for $H(C, \partial)$, the homology of $C$ with
respect to the operator $\partial$. Denote this basis for
$H(C, \partial)$ by ${\bf h}^w, \, w \in W$. Then its 
inverse image under the map $\psi_{\partial, S}$
is a basis for $\Ker (S)$. In other words, for each $w \in W$, set
\[
s^w \, = \, \psi_{\partial, S}^{-1}({\bf h}^w) \, \in \, 
\Ker (S).
\]
We will refer to the forms $s^w$, for $w \in W$, as the {\bf 
Kostant harmonic forms}.  The form $s^w$ is of pure degree
$2l(w)$, where $l(w)$ is the length of $w$.  In fact, it has bi-degree
$(l(w), l(w))$ in the bi-grading of $C$ defined 
by the complex structure $J$ on $K/T \cong G/B_{-}$, where 
$G$ is the complexification of $K$ and $B_{-}$
the Borel subgroup of $G$ corresponding to the 
Lie algebra $\fb_{-} = \fh + \fn_{-}$.

The rest of \cite{ko:63} is devoted to 
the proofs of several special properties of the forms $s^w, \, 
w \in W$. 
We list some of these properties here.

\begin{enumerate}
\item
The basis $\{[s^w]: \, w \in W \}$ of
$H^{\bullet}(K/T, \C)$ is, up to scalar multiples,
dual to the basis in the homology of 
$K/T$ defined by the Schubert varieties. In other words, if we use
$\Sigma_w$ to denote the Schubert cell in $K/T$ 
corresponding to $w$ (defined as the orbits of the
$B_{+}$-action on $K/T \cong G/B_{+}$,
where $B_{+}$ is the
Borel subgroup of $G$ corresponding to the Lie algebra
$\fb_{+} = \fh + \fn_{+}$), then
\[
\int_{\Sigma_{w_1}} s^w \, = \, 
\left\{ \begin{array}{ll} 0 & {\rm if} \, \, w_1 \neq w \\
\lambda_w  \neq 0& {\rm if} \, \, w_1 = w \end{array} \right..
\]
The number $\lambda_w$ was later calculated by
Kostant and Kumar  in \cite{k-k:integral}, where they also
treated the general Kac-Moody case.

\item
For $w \in W$, denote by $i_{w}: \Sigma_w \hookrightarrow
K/T$ the inclusion map.
Then, in fact, when 
$l(w) = l(w_1)$ but $w \neq w_1$, the form
$i_{w_1}^{*} s^w$ on $\Sigma_{w_1}$ is identically equal to zero.

\item
Kostant actually constructed harmonic forms on $G/P$, where
$P \subset G$ is any parabolic subgroup containing $B_{+}$. Denote by
$W_P \subset W$ the subgroup of $W$ corresponding to $P$ and
by $W^P$ the set of minimal length representatives
of the right coset space $W_P \backslash W$. Then for
each $w \in W^P$, Kostant constructed a $K$-invariant 
harmonic form $s^{w}_{P}$ on $G/P \cong K/(K \cap P)$.
Under the natural projection $\nu: G/B_{+} \rightarrow G/P$,
the harmonic form $s^{w}_{P}$ on $G/P$ goes to the harmonic form
 $s^w$ on $G/B_{+}$ for each $w \in W^P$.

\item
In \cite{k-k:integral}, Kostant and Kumar described 
 a connection between
the harmonic forms $s^w, w \in W$, and the
Schubert calculus on $K/T$. Namely, a certain $D$-matrix,
which encodes the structure constants of the ring structure on
$H^{\bullet}(K/T, \C)$ in the Schubert basis, can be constructed using 
these harmonic forms.
This is reviewed in Section \ref{sec_equi}.
\end{enumerate}

\bigskip
As we mentioned earlier, what we achieve in this
paper is a Poisson geometrical interpretation
of Kostant's operator $\partial$ and the forms $s^w$ for $w \in W$.  
This will allow us to give Poisson geometrical proofs of the above
properties of these forms.
We now explain our results and the organization of the paper in more detail.

\bigskip
The notion of Poisson harmonic forms is introduced in
Section \ref{sec_poi-harm}, where we also prove a few 
general facts about such forms that will be 
applied to the Kostant harmonic forms in later sections.

\bigskip
The definition of the Bruhat Poisson structure $\pinf$ on the flag
manifold $K/T$ is reviewed in Section \ref{sec_bruhat}.
The Koszul-Brylinski operator $\partial_{\pi_{\infty}, \mu_0} $
defined by $\pinf$ and a $K$-invariant volume form $\mu_0$ on
$K/T$ is
denoted by $\ponf$:
\[
\ponf \, : = \, i_{\pi_{\infty}} d \, - \, d i_{\pi_{\infty}}
\, + \, i_{\theta_0}.
\]
The modular vector field $\theta_0 := \theta_{\mu_0}$
in this case turns out to be 
the infinitesimal generator of the $T$-action on $K/T$ (by left
translations) in the direction of $iH_{\rho} \in \ft$
(see Section \ref{sec_bruhat}).
The operator $\ponf$ leaves invariant the subspace $C = 
\Omega(K/T, \C)^K$
of $K$-invariant differential forms on $K/T$, and on $C$, 
the two operators $d$ and $\ponf$ anti-commute. 

\bigskip
In Section \ref{sec_mixed}, using 
an identification given in (\ref{eq_ii-multi}):
\[
I_{\infty}: \, \, (\wedge^{\bullet}(\fn_{-} \oplus
\fn_{+}))^T \stackrel{\sim}{\lrw} C^{\bullet},
\]
the operator $\ponf$ is identified 
with the Chevalley-Eilenberg boundary
operator $b_{{\frak n}_{-} \oplus {\frak n}_{+}}$
for the Lie algebra $\fn_{-} \oplus \fn_{+}$:

\begin{itemize}
\item
Theorem \ref{thm_ponf-n-complex} {\it
As operators on $C = \Omega(K/T, \C)^K$, we have
\[
\ponf \, = \, I_{\infty} b_{{\frak n}_{-} \oplus {\frak n}_{+}}
I_{\infty}^{-1}: \, \, C^q \lrw C^{q-1}.
\]
}
\end{itemize}

\noindent
This fact is the key reason why the Bruhat Poisson structure 
should be related to Kostant's harmonic forms:
Recall that the Kostant operator $\partial$ is defined as
$\partial = -I b_{{\frak n}_{-} \oplus {\frak n}_{+}} I^{-1}$, where
$I$ is the identification in  (\ref{eq_i-intro})
between the same two spaces. The two identifications
$I$ and $I_{\infty}$ are different, but their
difference is expressed via the standard complex structure $J$
on $K/T \cong G/B_{-}$, namely
$I_{\infty} I^{-1} = J$ (Lemma \ref{lem_I-J}). The immediate
consequences of this important point are
the following three of our main theorems, where 
$C = \oplus_{p,q} C^{p,q}$ is the bi-grading on
$C$ defined by $J$:

\bigskip
\begin{itemize}
\item
{\bf Theorem \ref{thm_main-1}}
{\it 
As operators on $C = \Omega(K/T, \C)^K$, 
the Kostant operator $\partial$ and
the Koszul-Brylinski operator $\ponf$ are
related by
\[
\partial \, = \, J \ponf J^{-1}
\]
where $J$ is the standard complex structure on $K/T \cong
G/B_{-}$.
}

\item
{\bf Theorem \ref{thm_main-2}}  {\it 
A form  $\xi \in C^{p,q}$ is $\partial$-closed if and only if it is
$\ponf$-closed.
}

\item
{\bf Corollary \ref{cor_main-2}} {\it
A form $\xi \in C^{p,q}$ is
harmonic in the sense of Kostant if and only if it is Poisson
harmonic with respect to the Bruhat Poisson structure $\pinf$
(and any $K$-invariant volume form). 
In particular, the Kostant harmonic forms $s^w$, for $w \in W$,
are all Poisson harmonic.}
 \end{itemize}

\bigskip
Using this Poisson characterization of the forms $s^w$
for $w \in W$,
we immediately derive their first three properties
 listed earlier. This is done in Section \ref{sec_kos-thm}.

\bigskip
In Section \ref{sec_equi}, we consider the $S^1$-equivariant 
cohomology $H_{S^{1}}(K/T)$ of $K/T$, where the $S^1$-action 
on $K/T$ is defined by the element $iH_{\rho} \in \ft$, 
or, in other words, by the modular vector field $\theta_0$.
More precisely,
for each $w \in W$, we define the $\Omega(K/T, \C)$-valued
function on $\R$ by
\[
s^w(u) \, = \, i_{\exp_{\wedge} (-u \pi_{\infty})} s^w \, = \, 
s^w \, - \,u i_{\pi_{\infty}} s^w \, + \, 
{u^2 \over 2!} i_{\pi_{\infty} \wedge \pi_{\infty}} s^w \, + \, \cdots.
\]
Following from an easy fact about Poisson harmonic forms 
(Theorem \ref{thm_equi} in Section \ref{sec_poi-harm}),
we have

\bigskip
\begin{itemize}
\item
{\bf Theorem \ref{thm_main-3}} {\it 
Each $s^w(u)$ is $S^1$-equivariantly closed, and, 
up to scalar multiples, the set $\{[s^w(u)]: ~ w \in W \}$ is the 
$\C[u]$-Schubert basis 
for the $S^1$-equivariant cohomology of $K/T$.
}
\end{itemize}

\bigskip
This theorem is then used to give another proof of the 
connection between the Kostant harmonic forms and the 
Schubert calculus on $K/T$ that was found by Kostant and Kumar
in \cite{k-k:integral} 
(the fourth property of these forms listed earlier):
By evaluating the forms $s^w(u)$ at the $T$-fixed points on $K/T$,
we get a matrix $D$ that
encodes the structure constants of the $\C[u]$-algebra structure
on $\hus$ in the Schubert basis. Since $\hus$ specialized at
$u = 0$ is the ordinary de Rham cohomology $H^{\bullet}(K/T, \C)$,
we get an explicit description of the ring structure on 
$H^{\bullet}(K/T, \C)$ in terms of the Kostant harmonic forms. 
We remark that a closely related algorithm for computing the ring structure
was given in the seminal paper \cite{bgg:73}, and that the argument in
\cite{bgg:73} is based on consideration of certain $\C P^1$ bundles
over $G/B_{+}.$ We are not using $\C P^1$ bundles in any explicit way
to obtain our formulas, although we believe that recovering formulas
based on $\C P^1$ bundles from our arguments should not be very difficult.
We plan to pursue this topic in a future paper.

\bigskip
Section \ref{sec_family} is devoted to a 
family of Poisson (in fact symplectic) structures on 
$K/T$
denoted by $\pi_{\lambda}$, 
where the parameter $\lambda$ runs over all regular elements in
$\ft^*$, 
and the symplectic structure $\pi_{\lambda}$  comes from the dressing orbit
of $K$ in its dual group through the point $e^{-\lambda}$.
Very importantly, the Poisson structures $\pi_{\lambda}$ tend to
the Bruhat Poisson structure as $\lambda \rightarrow \infty$.
We study Poisson harmonic forms for $\pi_{\lambda}$. It turns
out that the Koszul-Brylinski operator
$\pal := \partial_{\pi_{\lambda}, \mu_0}$ for $\pi_{\lambda}$
 and the adjoint $d_{\ast, \lambda}$
of $d$ with respect to a certain Hermitian metric $h_{\lambda}$ are related by
$d_{\ast, \lambda} = J \pal J^{-1}$ (compare with
$\partial = J \ponf J^{-1})$.
Consequently, we have (Theorem \ref{thm_limit-S})
\[
 \partial \, = \, \lim_{t \rightarrow +\infty} 
d_{\ast,\lambda + t H_{\rho}} \hspace{.3in} {\rm and} \hspace{.3in}
S \, = \, \lim_{t \rightarrow +\infty} S_{\lambda + t H_{\rho}},
\]
where $S_{\lambda} = d d_{\ast,\lambda} + d_{\ast,\lambda} d$
is the usual (Hodge) Laplacian for $d$ with respect to the Hermitian
metric $h_{\lambda}$. Thus, although $\partial$ is in general not
the adjoint of $d$ with respect to any Hermitian metric, it is
the limit of a family of such operators.

The construction of the forms $s^w$ for $w \in W$ in \cite{ko:63}
relies on the fact that the two operators $d$ and $\partial$
are disjoint in the sense that
$\Im (d) \cap \Ker (\partial)
= \Im (\partial) \cap \Ker (d) = 0$.
Notice that  $d$ is always disjoint from 
to its adjoint operator with respect to any
Hermitian metric. Using Theorem \ref{thm_limit-S}
and some simple linear algebra arguments, we
give another proof of the disjointness of $d$ and $\partial$.
This is done in Section \ref{sec_another-proof}.

Finally, for each $\pi_{\lambda}$, we construct Poisson harmonic forms
$s^{w}_{\lambda}, w \in W$.
We show that $s^{w}_{\lambda + t H_{\rho}} \rightarrow s^w$
as $t \rightarrow +\infty$.  Thus the Kostant harmonic forms
are limits of usual Hodge harmonic forms.

\bigskip
We conclude this introduction by mentioning
a few other works that are related to this paper.

\bigskip
The notion of the modular vector field $\theta_{\mu}$
of a Poisson structure $\pi$ associated to a volume form $\mu$,
which plays a central role in our work
(see Sections \ref{sec_poi-harm}
and \ref{sec_equi}), has been actively investigated only recently:
although it appeared in some earlier work of Koszul, it was rediscovered
by Weinstein \cite{we:modular} and independently by Brylinski and Zuckerman
\cite{b-z:modular} in 1995. See also  
the paper \cite{po:poi} by Polishchuk.
A related paper is \cite{e-l-w:modular}, in 
which we study modular vector fields for arbitrary Lie algebroids. 
Results in \cite{e-l-w:modular} show that there is a notion of Poisson 
harmonic forms on a Poisson manifold $P$
for each rank $1$ representation of the cotangent bundle Lie 
algebroid $T^*P$. The notion of Poisson harmonic forms we use in this
paper corresponds to the trivial rank $1$ representation of $T^*P$, while
that in \cite{by:poi} corresponds to the representation of $T^*P$
on the canonical line bundle of $P$. See \cite{e-l-w:modular} 
\cite{po:poi} and \cite{xu:modular} for more details.

\bigskip
The Bruhat Poisson structure $\pinf$ on $K/T$ is an example
of a $(K, \pi)$-homogeneous Poisson structure, where
$\pi$ is the Poisson structure on $K$ given in (\ref{eq_pi-on-K})
which makes $(K, \pi)$ into a Poisson Lie group.
Poisson homogeneous spaces for Poisson Lie groups 
were first studied by Drinfeld \cite{dr:homog}, where he shows that
a certain Lie algebra determines the Poisson structure on the underlying
homogeneous space. For the Bruhat Poisson structure, this Lie
algebra is $\ft + \fn_{+}$ as a real Lie subalgebra of $\fg$. 
This can be regarded as the first hint why Kostant's work in
\cite{ko:63} and the Bruhat Poisson structure $\pinf$ 
should be related: the Lie algebra $\fn_{+}$ appears in both
settings. In \cite{lu:homog}, we extend the 
work of Drinfeld and show that the (real) invariant Poisson cohomology 
of the Poisson homogeneous space is isomorphic to 
certain relative Lie algebra cohomology of the Lie algebra defined
by Drinfeld--it is
$ H^{\bullet}(\ft + \fn_{+}; \ft) 
\cong (H^{\bullet} (\fn_{+}))^T$
in the case of $\pinf$ (note that $\fn_{+}$
is considered as a real Lie algebra here).
It is clear at this point that Kostant's results in 
\cite{ko:63} are statements about the Bruhat Poisson
structure and its (complex) $K$-invariant Poisson cohomology. 
This is the starting point
of our project to give Poisson interpretations of
Kostant's work.

\bigskip
The work in \cite{lu:coor} is our first attempt
to relate Kostant's harmonic forms to the geometry
of the Bruhat Poisson structure.
In \cite{lu:coor}, for each $w \in W$,
we explicitly express the Kostant harmonic form $s^w$ 
as a Duistermaat-Heckman type volume form on the cell
$\Sigma_w$ 
(as a symplectic leaf of $\pinf$).
In fact, the form $s^w$, as well as some other related
quantities such as the moment map for the $T$ action on $\Sigma_w$, 
are written
down in \cite{lu:coor} by explicit  formulas in
certain Bott-Samelson type coordinates on the cell. The calculation there
for the integral $\lambda_w = \int_{\Sigma_w} s^w$
takes only a couple of lines by using the explicit formula in coordinates.
Some of the results in \cite{lu:coor}, for example,
Theorem 4.3 there, follow from the global Poisson properties
of the forms $s^w$ we find in this paper. This is
explained in Section \ref{sec_kos-thm}.

\bigskip 
Ideas in Section \ref{sec_family} come from 
\cite{lu:cdyb}, where
a family of Poisson structures $\pi_{\tx, \lambda}$ on $K/T$,
of which the $\pi_{\lambda}$'s considered in Section \ref{sec_family}
constitute only a subfamily,
is studied. The Poisson structures in this family all come
from the solutions to the Classical Dynamical Yang Baxter Equation
for the pair $(\fg, \fh)$ that have recently been classified
by Etingof and Varchenko \cite{e-v:cdyb}, and this family exhausts
all $(K, \pi)$-homogeneous Poisson structures on $K/T$.
The fact that 
the Bruhat Poisson structure is the limit of this family as $\lambda
\rightarrow \infty$ is first proved in \cite{lu:cdyb}. 
Many Poisson geometrical properties of 
members of the family $\pi_{\tx, \lambda}$, such as their symplectic 
leaves, their  $K$-invariant Poisson cohomology, and the Drinfeld
Lie algebras are studied in \cite{lu:cdyb}.

\bigskip
Papers \cite{lu:homog} \cite{lu:coor} \cite{lu:cdyb} and this one
can be regarded as a series, where the Poisson geometrical 
properties of certain Poisson structures are studied, while 
applications to Lie theory go along with the study.
More papers are planned, where connections between
the Bruhat Poisson structure $\pinf$ and the $T$-equivariant cohomology
(instead of the $S^1$-equivariant cohomology) of $K/T$, connections
between the Bruhat Poisson structures and the BGG operators, Poisson 
geometrical aspects of the Bott-Samelson resolutions for
Schubert varieties, and 
the Poisson (co)homology (as opposed to 
the $K$-invariant Poisson cohomology) of the family 
$\pi_{\tx, \lambda}$ will be studied.

\bigskip
{\bf Acknowledgement} We would like to thank Bert Kostant for
explaining to us his results in \cite{ko:63} and for 
his constant encouragement for this work. We also thank
Alan Weinstein and Dale Peterson for useful discussions.
The second author is grateful to the Hong Kong University
of Science and Technology for its hospitality.

\section{Poisson harmonic forms}
\label{sec_poi-harm}

Recall that a Poisson structure on a manifold $P$ is
a bi-vector field $\pi$ on $P$ satisfying
\begin{equation}
\label{eq_pi}
[\pi, ~ \pi] ~ = ~ 0,
\end{equation}
where $[\pi, ~ \pi]$ is the Schouten bracket of $\pi$
with itself. The Schouten bracket is reviewed in the Appendix.

\bigskip
Denote by $\chi^q(P)$ the space of $q$-vector fields
on $P$, i.e., smooth sections of the vector bundle 
$\wedge^q TP$, and by $\Omega^q(P)$ the space of $q$-differential
forms on $P$. Define
\beqa
\delta_{\pi} & = & [\pi, \, \bullet]: \, \,
\chi^q(P) \lrw \chi^{q+1}(P),\\
\partial_{\pi} & = & i_{\pi} d \, - \, d i_{\pi}: \, \, 
\Omega^{q}(P) \lrw \Omega^{q-1}(P).
\eeqa
It follows from the graded Jacobi identity for 
the Schouten bracket that
$\delta_{\pi}^{2} = 0$, and it follows from 
(\ref{eq_p-bra}) in the Appendix that $\partial_{\pi}^{2} = 0$.
The {\bf Poisson
cohomology} \cite{li:poi} of $(P, \pi)$ is defined 
to be the cohomology of the cochain complex
$(\chi^{\bullet}(P), \delta_{\pi})$ and is denoted by
$H^{\bullet}(P, \pi)$. The {\bf Poisson homology} 
\cite{by:poi} of $(P, \pi)$
is defined to be the homology of the chain complex
$(\Omega^{\bullet}(P), \partial_{\pi})$.
See \cite{vs:poi-book} for more details on Poisson cohomology.

\bigskip
Assume now that $(P, \pi)$ is an orientable Poisson manifold
of dimension $n$ and that
$\mu$ is  a volume form on $P$. The map
\[
\backl \, \mu: \, \chi^q(P) \lrw \Omega^{n-q}(P): \, \,
X \Map X \backl \, \mu \, = \, i_X \mu
\]
is an isomorphism of vector spaces, 
where $X \backl \mu = i_X \mu$ 
is the contraction of $X$ with $\mu$. \footnote{Throughout this paper,
we will use both $X \backl \xi$ and $i_X \xi$ to denote
the contraction of a multi-vector field $X$ with a differential
form $\alpha$.}

\begin{lem}
\label{lem_delta-partial}
For any $q$-vector field $X$ on $P$, we have
\[
(\delta_{\pi} X) \backl \, \mu \, = \, (-1)^{q-1}
(\partial_{\pi} + i_{\theta{\mu}}) (X \backl \, \mu),
\]
where $\theta_{\mu}$ is the unique vector field on $P$ 
such that
\[
\theta_{\mu} \backl \, \mu \, = \, d(\pi \backl \, \mu) \, = \, 
-\partial_{\pi} \mu.
\]
\end{lem}

\bigskip
\noindent
{\bf Proof.} 
It follows from (\ref{eq_i-bra}) that
\beqa
(\delta_{\pi} X) \backl \, \mu & = & i_{[\pi, \, X]} \mu \\
& = & (-1)^{q-1}(\partial_{\pi} i_{X} \mu
\, - \, (-1)^{q}i_X \partial_{\pi} \mu) \\
& = & (-1)^{q-1} (\partial_{\pi} i_{X} \mu \, + \,
(-1)^{q} i_X i_{\theta_{\mu}} \mu) \\
& = & (-1)^{q-1} (\partial_{\pi} \, + \, i_{\theta_{\mu}} )
(X \backl \mu).
\eeqa
\qed

\begin{dfn} \cite{we:modular}
\label{dfn_modular-we}
{\em
The vector field $\theta_{\mu}$ is called the
{\bf modular vector field} of $\pi$ associated to
the volume form $\mu$.}
\end{dfn}
 
The modular vector field is always a Poisson vector field, i.e.,
$[\theta_{\mu}, \, \pi] = 0$. The class it defines in $H^{1}(P, \pi)$
is independent of the choice of the volume form $\mu$, and
this class is called the modular class of $\pi$ ( 
\cite{we:modular} 
 \cite{e-l-w:modular} \cite{xu:modular}).
 
\begin{nota}\label{nota:partial}
{\em
We set
\[
\partial_{\pi, \mu} \, = \, \partial_{\pi} \, + \, i_{\theta_{\mu}}
\, = \, i_{\pi} d \, - \, d i_{\pi} \, + \, i_{\theta_{\mu}}:
\, \, \Omega^q(P) \lrw \Omega^{q-1}(P).
\]
}
\end{nota}

\bigskip
\begin{dfn}
\label{dfn_poi-harm}
{\em
We say that a differential form $\xi$ on $P$ is {\bf Poisson harmonic
with respect to the Poisson structure $\pi$ and the
volume form $\mu$} if it satisfies
\[
d \xi \, = \, 0 \hspace{.5in} {\rm and} \hspace{.5in}
 \partial_{\pi, \mu} \xi \, = \, 
0.
\]
We will also say that such forms are $(d, \partial_{\pi, \mu})$-harmonic.
}
\end{dfn}

\begin{exam}
\label{exam_1}
{\em
The volume form $\mu$ is always $(d, \partial_{\pi, \mu})$-harmonic,
and a general top degree form $f \mu$ is $(d, \partial_{\pi, \mu})$-harmonic
if and only if $f$ is a Casimir function on $P$, i.e., 
$df \backl \pi = 0$. This is easily seen from Lemma 
\ref{lem_delta-partial}. 
}
\end{exam}

\begin{exam}
\label{exam_2}
{\em
Let $(P, \omega)$ be a symplectic manifold. Then there is a unique 
non-degenerate Poisson bi-vector field $\pi$ on $P$ characterized by
\begin{equation}
\label{eq_omega}
\omega( \tilde{\pi} (\alpha), \, \tilde{\pi}(\beta)) \, = \, \pi
(\alpha, \, \beta)
\end{equation}
for any $1$-forms $\alpha$ and $\beta$, where $\tilde{\pi} (\alpha)
= \alpha \backl \pi$. Let $\mu_{\omega} = {\omega^m \over m!}$,
where $2m = \dim P$,
be the Liouville volume form on $P$ defined by $\omega$. Then
the modular vector field of $\pi$ associated to
$\mu_{\omega}$ is zero. In general, let
$\mu = f \mu_{\omega}$ be an arbitrary volume form
on $P$, where $f$ is a nowhere vanishing function on $P$.
Then the modular vector field of $\pi$ associated to $\mu$ is 
$\theta_{\mu} = -d(\log |f|) \backl \pi$, or, in other words,
the Hamiltonian vector field of the function $\log |f|$.
As a particular case of Example \ref{exam_1}, when $P$ is connected, 
a top degree form $s$ on $P$ is $(d, \partial_{\pi, \mu})$-harmonic
if and only if $s = c \mu$ for a constant function $c$.
This simple fact will be used in Section \ref{sec_kos-thm}.
}
\end{exam}

\begin{rem}
\label{rem_k-b-operator}
{\em
Note that the volume form $\mu$ enters into our definition of
Poisson harmonic forms. 
The operator $\partial_{\pi} = i_\pi d - d i_{\pi}$ was first
studied by Koszul \cite{kz:schouten} and Brylinski \cite{by:poi}
and is generally called the Koszul-Brylinski operator
associated to $\pi$. It is thus appropriate to call
$ \partial_{\pi, \mu} = \partial_{\pi} + i_{\theta_{\mu}}$
the {\bf Koszul-Brylinski operator associated to $\pi$ and $\mu$}. 
Of course, 
it also satisfies $\partial_{\pi, \mu}^{2} = 
0$, but its homology is isomorphic to the Poisson cohomology of $\pi$. 
Poisson harmonic forms are defined in \cite{by:poi} to be 
$(d, \partial_{\pi})$-harmonic. One can show by using our results
in Section \ref{sec_kos-thm} that
for the Bruhat Poisson structure on 
the flag manifold $K/T$, a non-zero de Rham cohomology class 
has a  $(d, \partial_{\pi})$-harmonic 
representative only when it is in degree $0$
while this is not the case for 
$(d, \partial_{\pi, \mu})$-harmonic forms. 
The following question is a modification
of a question asked by Brylinski in \cite{by:poi}.
 }
\end{rem}

\begin{ques}
\label{ques_harm}
Given an orientable manifold $P$ with
a Poisson structure $\pi$ and a volume form $\mu$,
 does every de Rham cohomology class of $P$ have a
representative that is Poisson harmonic
with respect to $\pi$ and $\mu$?
\end{ques}

\bigskip
When $\pi$ is non-degenerate and when 
$\mu$ is the Liouville
volume form corresponding to $\omega$ given by (\ref{eq_omega}), 
the modular vector field
$\theta_{\mu}$ is zero, and thus 
$\partial_{\pi, \mu} = i_{\pi} d - d i_{\pi}$
is the same as $\partial_{\pi}$. A result of
O. Mathieu \cite{mt:harm} says that the answer to 
Question \ref{ques_harm} is ``yes" if and only if
for any $k \leq m = {1 \over 2} \dim P$, the
cup product by $[\omega]^k: H^{m-k}(P) \rightarrow
H^{m+k}(P)$
is an isomorphism. A simpler proof of this fact is given in \cite{yd:harm}.
 
Our results in Section \ref{sec_kos-thm} will show that for the
Bruhat Poisson structure on the flag manifold $K/T$
and a $K$-invariant volume form $\mu_0$ on $K/T$, the
answer to Question \ref{ques_harm} is ``yes". 
Other examples for which the answer
to Question \ref{ques_harm} is ``yes" are given in 
Section \ref{sec_poi-hodge}. See Remark \ref{rem_rep-ok-pl}.

\bigskip
We now look at the connection between Poisson harmonic forms and
$S^1$-equivariant cohomology.

\bigskip
Let again $(P, \pi)$ be a Poisson
manifold. Introduce the elements
\begin{eqnarray}
\label{eq_exp-pi}
\exp_{\wedge} \pi ~ = ~ 1 + \pi + {\frac{1}{2!}}
\pi \wedge \pi + {\frac{1}{3!}} \pi \wedge \pi \wedge \pi +
\cdots \in \chi (P)\\
\label{eq_exp}
\exp_{\wedge} (-\pi) ~ = ~ 1 - \pi + {\frac{1}{2!}}
\pi \wedge \pi - {\frac{1}{3!}} \pi \wedge \pi \wedge \pi +
\cdots \in \chi(P).
\end{eqnarray}
They are inverse to each other with respect to the
wedge product on $\chi(P)$. Correspondingly, the operator
\[
i_{\exp_{\wedge} \pi}: ~ \Omega^k(P) \lrw
\Omega^k(P) ~ \oplus ~ \Omega^{k-2}(P) ~ \oplus ~
\Omega^{k-4}(P) ~ \oplus \cdots
\]
given by
\begin{equation}
\label{eq_iep}
\iep (\alpha) ~ = ~ \alpha ~ + ~ i_{\pi} \alpha
~ + ~ {\frac{1}{2!}} i_{\pi \wedge \pi} \alpha
~ + ~
{\frac{1}{3!}} i_{\pi \wedge \pi \wedge \pi} \alpha
~ + ~ \cdots
\end{equation}
has the operator $\ienp$ as its inverse. Since
$i_{\pi \wedge \pi} = i_{\pi}^2$, we have,
\[
\iep ~ = ~ {\rm Exp} (i_{\pi}),
\]
where ${\rm Exp} (i_{\pi})$ is the exponential of the operator
$i_{\pi}$. It is the sum of only finitely many terms.

The following proposition is also proved by Polishchuk \cite{po:poi}.

\begin{prop}
\label{prop_iep}
On $\Omega(P)$, we have
\begin{equation}
\label{eq_d-pdp}
d ~ + ~ \pdp ~ = ~ (\iep) ~ d ~ (\iep)^{-1}.
\end{equation}
\end{prop}
 
\noindent
{\bf Proof.}
Denote by $[ ~~~]$ the commutator bracket between operators on
$\Omega(P)$. We know from (\ref{eq_i-bra}) that
\[
[i_{\pi}, ~ \pdp] ~ = ~  [i_\pi, ~ [i_\pi, ~ d]] ~ = ~ 0.
\]
Thus
\beqa
\iep ~ d ~ \ienp & = & {\rm Exp} (i_\pi) \, d \,{\rm Exp}
(-i_\pi)\\
& = & d ~+ ~  [i_{\pi}, ~ d] ~ + ~
{\frac{1}{2!}} [i_\pi, ~ [i_\pi, ~ d]]
~ + ~ {\frac{1}{3!}}[i_\pi, ~ [i_\pi, ~ [i_\pi, ~ d]]] ~ + ~
\cdots\\
& = &  d ~ + ~ \pdp.
\eeqa
\qed

Assume now that $P$ is orientable. Let $\mu$ be a volume form,
and let $\theta_{\mu}$ be the modular vector field
of $\pi$ associated to $\mu$.
Recall that a $\theta_{\mu}$-equivariantly closed differential
form on $P$ is an $\Omega(P)$-valued polynomial
function $\xi(u)$ on $\R$ such that
\[
d \xi(u) \, + \, i_{u \theta_{\mu}} \xi(u) \, = \, 0, \hspace{.5in}
\forall u \in \R.
\]
When the vector field $\theta_{\mu}$
integrates to an $S^1$-action on $P$, such forms are said to be
$S^1$-equivariantly closed.
The following observation is the key to our construction of
$S^1$-equivariantly closed forms on the flag manifolds in Section 
\ref{sec_equi}.

\begin{thm}
\label{thm_equi}
For any differential form  $\xi$ on $P$ of pure degree, define
\[
\xi(u) \, = \, i_{\exp_{\wedge} (-u \pi)} \xi, 
\]
for $u \in \R$. Then $\xi$ is $(d, \partial_{\pi, \mu})$-harmonic if
and only if $\xi(u)$ is $\theta_{\mu}$-equivariantly closed.
\end{thm}

\bigskip
\noindent
{\bf Proof.} By Proposition \ref{prop_iep}, we have
\beqa
 d \xi(u) & = & d i_{\exp_{\wedge} (-u \pi)} \xi \\
& = & i_{\exp_{\wedge} (-u \pi)}  (d \, + \, 
\partial_{u \pi}) (\xi) \\
& = & i_{\exp_{\wedge} (-u \pi)} (d \, + \, u \partial_{\pi, \mu}
\, - \, u i_{\theta_{\mu}})(\xi).
\eeqa
Thus 
\[
d \xi(u) \, + \, i_{u \theta_{\mu}} \xi(u) \, = \
i_{\exp_{\wedge} (-u \pi)} (d \, + \, u \partial_{\pi, \mu})(\xi).
\]
Since $\xi$ is of pure degree, we see that 
$ d \xi(u) \, + \, i_{u \theta_{\mu}} \xi(u) = 0$ for all $u \in \R$
if and only if $d \xi = 0$ and $\partial_{\pi, \mu} \xi = 0$, i.e.,
if and only if $\xi$ is $(d, \partial_{\pi, \mu})$-harmonic.
\qed

\begin{exam}
\label{exam_3}
{\em
Let $\pi$ be a non-degenerate Poisson structure on a manifold $P$
of dimension $2m$. Assume that $\theta$ is a Hamiltonian vector field
on $P$ with a Hamiltonian function $\phi$, i.e., $\theta = 
-d \phi \backl \pi$. Consider the volume form
$ \mu = e^{\phi} {\omega^n \over n!}$, where $\omega$ is the
symplectic form on $P$ defined by $\pi$ as given in (\ref{eq_omega}).
Then the modular vector field of $\pi$ associated to $\mu$
is $\theta$, and thus by Example \ref{exam_1} and Theorem 
\ref{thm_equi}, we have the $\theta$-equivariantly closed
form $\mu(u)$ on $P$ given by
\beqa
\mu (u) &  = & i_{\exp_{\wedge}(-u \pi)} \mu \\
& = & e^{\phi} \left({\omega^m \over m!} \, - \, 
u {\omega^{m-1} \over (m-1)!} \, + \, u^2 {\omega^{m-2} \over (m-2)!}
\, + \cdots + \,
(-1)^{m-1} u^{m-1} \omega \, + \, (-1)^m u^m \right) \\
& = & (-u)^m \exp(\phi - {\omega \over u}).
\eeqa
This is the well-known $\theta$-equivariant
form used by Atiyah and Bott in \cite{a-b:equi} to
prove the Duistermaat-Heckman integral formula from 
the localization theorem for equivariant cohomology.
}
\end{exam}

\section{The Bruhat Poisson structure and Kostant's harmonic forms}
\label{sec_bruhat-kostant}

\subsection{The Bruhat Poisson structure $\pinf$ and the Koszul-Brylinski
operator $\ponf$}
\label{sec_bruhat}

Let $K$ be a compact semi-simple Lie group and $T \subset K$
a maximal torus. Denote by $\fk$ and $\ft$ the Lie algebras of 
$K$ and $T$ respectively.
Let $\fg$
be the complexification of $\fk$. Then the complexification
$\fh$ of $\ft$ is a Cartan subalgebra of $\fg$. 
Denote by $\Phi$ the set of all roots
of $\fg$ with respect to $\fh$ and by $\Phi^{+}$ a choice of
positive roots. We will also write $\alpha > 0$ for
$\alpha \in \Phi^+$. 

Let $\ll ~ \gg$ be the Killing form of $\fg$. For each
positive root $\alpha$, denote by $H_{\alpha}$ the image of $\alpha$
under the isomorphism $\fh^* \rightarrow \fh$ via $\ll ~ , ~ \gg$,
i.e., for any $H \in \fh$,
\[
\ll H_{\alpha}, ~ H \gg  ~ = ~ \alpha(H).
\]
Choose $\Xa, \Ya \in \fk$ such that
$\ll \ea, ~ \eb \gg 
~ = ~ 1$, where
\[
\ea \, = \, {1 \over 2}(\Xa \, - \, i\Ya), \hspace{.3in} 
{\rm and} \hspace{.3in} \eb \, = \, -{1 \over 2}
(\Xa \, + \, i \Ya)
\]
are root vectors for $\alpha$ and $-\a$ respectively.
It follows that $[\ea, \eb ] = H_{\alpha}$.
Set 
\begin{equation}
\label{eq_r}
r \,\, = \,\, {1 \over 4} \sum_{\alpha > 0} \Xa \wedge \Ya \, \, 
\in \,\, \fk \wedge \fk.
\end{equation}
This is the well-known $r$-matrix for the Lie algebra $\fk$.
Let $r^{R}$ and $r^{L}$ be respectively the 
right and left invariant bi-vector fields on
$K$ with values $r$ at the identity element $e$.
Define a bi-vector field $\pi$ on $K$ by
\begin{equation}
\label{eq_pi-on-K}
\pi \, = \, r^{R} \, - \, r^{L}.
\end{equation}
Then $\pi$ is a Poisson structure on $K$, and $(K, \pi)$
becomes a Poisson Lie group \cite{lu-we:poi}.

\bigskip
Since the Poisson structure $\pi$ on $K$ is invariant by the
right translations by elements in $T$, there is a unique
Poisson structure on $K/T$, which will be 
denoted by $\pi_{\infty}$ (the reason for
this notation will become clear in Section \ref{sec_family}), 
such that the natural projection
$(K, \pi) \rightarrow (K/T, \pi_{\infty})$
is a Poisson map. This Poisson structure
$\pi_{\infty}$ is called the Bruhat Poisson structure on $K/T$.
It is so named because its
symplectic leaves are
exactly all the Bruhat cells in $K/T$ 
\cite{lu-we:poi}. In \cite{lu:coor}, explicit
formulas for $\pinf$ on each Bruhat cell are written down in coordinates.

\bigskip
Note that although the Bruhat Poisson structure is
invariant under the left translations by 
elements in $T$, it is not invariant for elements in $K$. In
fact, the action map
\begin{equation}
\label{eq_action}
K \, \times K/T \lrw K/T: \, \, (k_1, \, kT) \Map k_1 kT
\end{equation}
is a Poisson map, where $K \times K/T$ is equipped with the
direct product Poisson structure $\pi \oplus \pinf$. This
is an example of a $(K, \pi)$-homogeneous Poisson space
\cite{dr:homog}.
The classification of all $(K, \pi)$-homogeneous
Poisson spaces is given by Karolinsky \cite{ka:homog-compact}.
Some geometric properties of all $(K, \pi)$-homogeneous
Poisson structures on $K/T$ are studied in \cite{lu:cdyb}.
A family of $(K, \pi)$-homogeneous Poisson structures
$\pi_{\lambda}$ on $K/T$ will be discussed in Section \ref{sec_family}.

\bigskip
Assume now that $\mu_0$ is a $K$-invariant volume form
on $K/T$. 
The modular vector field $\theta_{\mu_0}$
of $\pinf$ associated to $\mu_0$ is shown
in \cite{e-l-w:modular} to be the
infinitesimal generator of the $T$ action on
$K/T$ in the direction of the element $iH_{\rho}
\in \ft$, where $2 \rho = \sum_{\a > 0} \a \in i \ft$.

\bigskip
\begin{nota}
\label{nota_theta0}
{\em
Throughout this paper, the modular vector field
$\theta_{\mu_0}$ 
of $\pinf$ associated to $\mu_0$ will be denoted 
by $\theta_0$. Thus
\begin{equation}
\label{eq_theta-0}
\theta_0 (kT) \, = \, {d \over dt}|_{t=0} \exp(t iH_{\rho})k T.
\end{equation}
Notice that $\theta_0$ is independent of the choice of a
$K$-invariant volume form.
}
\end{nota}

\bigskip
Consider now the Koszul-Brylinski operator defined by $\pinf$ and $\mu_0$
\[
\partial_{\pi_{\infty}, \mu_0} \, = \, i_{\pi_{\infty}} d \, - \,
d i_{\pi_{\infty}} \, + \, i_{\theta_0}: \, \, \Omega^q(K/T) \lrw 
\Omega^{q-1}(K/T).
\]

\begin{nota}
\label{nota_partial-infty}
{\em
For notational simplicity, we set
\[
\partial_{\infty} \, = \, \partial_{\pi_{\infty}, \mu_0}: \, \,
\Omega^q(K/T) \lrw \Omega^{q-1}(K/T).
\]
}
\end{nota}
 
\bigskip
By Lemma \ref{lem_delta-partial}, the operator $\ponf$
is related to the operator
$\delta_{\pi_{\infty}} = [\pinf, \bullet]$ on $\chi(K/T)$
via
\begin{equation}
\label{eq_pinf-ponf}
(\delta_{\pi_{\infty}} X) \backl \, \mu_0 \, = \, (-1)^{q-1} \ponf
(X \backl \, \mu_0)
\end{equation}
for
$X \in \chi^q(K/T)$. Clearly,
\begin{equation}
\label{eq_d-ponf}
d \ponf \, + \ponf d \, = \, L_{\theta_0},
\end{equation}
where $L_{\theta_0}$ is the Lie derivative operator by the vector
field $\theta_0$.
The homology of the chain complex $(\Omega^{\bullet}(K/T), \, \ponf)$,
which is isomorphic to the Poisson cohomology of $\pinf$ by 
(\ref{eq_pinf-ponf}), can be shown to be isomorphic to the Lie algebra 
$\fn$-cohomology with coefficients in a certain infinite
dimensional module. 
The calculation of this cohomology will be carried
out in a separate paper. For now, we restrict ourselves
to the space of $K$-invariant differential forms on $K/T$.

\subsection{The mixed complex $(\Omega(K/T)^K, \, d, \, \ponf)$}
\label{sec_mixed}

\bigskip
Denote by $\chi(K/T)^K$
the space of $K$ invariant (real) multi-vector fields.
A general fact about Poisson actions (see, for example,
 $\S 7$ in  \cite{lu:homog}) is that the subspace
$\chi(K/T)^K \subset \chi(K/T)$ is invariant under
the operator $\delta_{\pinf} = [\pinf, \bullet]$. The
cohomology of the cochain complex
$(\chi(K/T)^K, \, \delta_{\pinf})$ is called the
{\bf $K$-invariant Poisson cohomology} of $(K/T, \, \pinf)$
(see \cite{lu:homog}).

Denote by $\Omega(K/T)^K$ the space of (real)
$K$-invariant  differential forms on
$K/T$. Then since $\mu_0$ is $K$-invariant, the map
\[
\backl \, \mu_0: \, \, \chi^q(K/T)^K \lrw \Omega^{n-q}(K/T)^K:
\, \, X \Map X \backl \, \mu_0
\]
is an isomorphism of vector spaces. We know from 
(\ref{eq_pinf-ponf}) that the operator
$\ponf$
on $\Omega(K/T)$ leaves the subspace $\Omega(K/T)^K$
invariant. It follows from (\ref{eq_d-ponf}) that
\[
d \ponf \, + \, \ponf d \, = \, 0
\]
as operators on $\Omega(K/T)^K$.
Thus on the graded vector space
$\Omega(K/T)^K  =  \oplus_q \Omega^q(K/T)^K$, we have
two anti-commuting  operators $d$ and $\ponf$, 
of degrees $1$ and $-1$ respectively,
such that $d^2 = 0$ and $\partial_{\infty}^{2} = 0$. 
In other words, we have a mixed complex
$(\Omega(K/T)^K, \, d, \, \ponf)$ \cite{lo:cyclic}.

\bigskip
In Section \ref{sec_partial}, we will compare the operator
$\ponf$ with the operator $\partial$ introduced by Kostant
in \cite{ko:63}. Since Kostant's operator $\partial$ 
is the Chevalley-Eilenberg boundary 
operator for a certain Lie algebra, we will need to
first identify $\ponf$ as a such. Since Kostant
considered complex valued differential forms on $K/T$, we
will also need to complexify the mixed complex
$(\Omega(K/T)^K, \, d, \, \ponf)$. These two tasks
occupy this section.

\bigskip
We first recall the Manin triple corresponding to the Poisson
Lie group $(K, \pi)$.

\bigskip
Consider the real Lie subalgebra $\fa + \fn$ of $\fg$, where
$\fa  =  i \ft$ and  $\fn = \fn_{+}$ is the subalgebra of $\fg$
spanned by all the positive root vectors (but
considered as a real Lie algebra here). 
Then
\[
\fg \, = \, \fk \, + \, \fa \, + \, \fn
\]
is an Iwasawa decomposition of $\fg$ (as a real semi-simple
Lie algebra). The triple $(\fg, \fk, \fa + \fn)$,
together with twice the imaginary part of the Killing form
$\la \, , \, \ra := 2 {\rm Im} \ll \, , \, \gg$, is the Manin triple
for the Poisson Lie group $(K, \pi)$ \cite{lu-we:poi}.

\bigskip
Set $\eo = eT \in K/T$, where $e$ is the identity element in $K$.
Identify
\[
\te \, \cong \, \fk / \ft.
\]
The pairing $2{\rm Im} \ll \, , \, \gg$ between $\fk$ and $\fa + \fn$ 
then induces a non-degenerate pairing between $\fk / \ft$ and $\fn$,
and we get an identification (of real vector spaces)
\begin{equation}
\label{eq_ii}
I_{\infty}: \, \, \fn \, \stackrel{\sim}{\lrw} \, \tse: \, \,
(\ii(x), \, y) \, = \, 2 {\rm Im} \ll x, \, y \gg
\end{equation}
for $x \in \fn$ and $y \in \fk / \ft \cong \te$. 
Note that $I_\infty$ is $T$-equivariant,
 where $T$ acts on $\fn$
by the Adjoint action and on $\tse$ by the linearization at
$\eo$ of the $T$ action on $K/T$. Since
\[
\Omega(K/T)^K \, \cong \, \wedge(\tse)^T,
\]
we get an identification, still denoted by $\ii$:
\[
\ii: \, \, (\wedge^q \fn)^T \, \stackrel{\sim}{\lrw} \,
\Omega^q(K/T)^K
\]
for each $0 \leq q \leq n = \dim_{\Bbb R} \fn = \dim_{\Bbb R}(K/T)$.
Let
\[
\iis: \, \, \te \lrw \fn^*
\]
be the dual of $\ii$. It gives rise to
\[
\iis: \, \, \chi^q(K/T)^K\cong (\wedge^q \te )^T
\stackrel{\sim}{\lrw} \, 
(\wedge^q \fn^*)^T.
\]
Let
\[
d_{\frak n}: \, \, \wedge^q \fn^* \, \lrw \, \wedge^{q+1} \fn^*
\]
be the Chevalley-Eilenberg coboundary operator 
for the Lie algebra $\fn$ (considered as a real Lie algebra). Its
restriction to the subspace $(\wedge \fn^*)^T \subset \wedge \fn^*$
is also denoted by $d_{\frak n}$. Let
\[
b_{\frak n}: \, \, \wedge^q \fn \, \lrw \, \wedge^{q-1} \fn: \, \, 
b_{\frak n} (x_1 \wedge x_2 \wedge \cdots \wedge x_q) 
\, = \, \sum_{i<j} (-1)^{i+j} [x_i, x_j] \wedge x_1 \wedge \cdots
\wedge \hat{x}_i \wedge \cdots \wedge \hat{x}_j \wedge \cdots \wedge
x_q
\]
be the Chevalley-Eilenberg boundary operator, 
and its restriction to $(\wedge \fn)^T \subset \wedge \fn$ will also
be denoted by $b_{\frak n}$. The two operators $d_{\frak n}$ and
$b_{\frak n}$ are dual to each other. We remark that 
our definition of the Chevalley-Eilenberg boundary operator
for a Lie algebra differs from that used by Kostant in \cite{ko:63}
by a minus sign.

\bigskip
\begin{thm}
\label{thm_ponf-n}
We have
\begin{equation}\label{eq_ponf-d}
(\iis)^{-1}  d_{\frak n} \iis \, = \, \delta_{\pi_{\infty}}
\, = \, [\pinf, \bullet]
\end{equation}
as degree $1$ operators on $\chi(K/T)^K$, and
\begin{equation}
\label{eq_ponf-partial}
\ii b_{\frak n} I_{\infty}^{-1} \, = \, \ponf
\end{equation}
as degree $-1$ operators on $\Omega(K/T)^K$.
\end{thm}

\bigskip
A proof of a general fact about $(K, \pi)$-homogeneous
Poisson structures on $K/T$, of which Theorem \ref{thm_ponf-n}
is a special case, is given in \cite{lu:cdyb}. Here, we give a direct 
proof for the sake of completeness.

\bigskip
We first prove a lemma about the Poisson Lie group $(K, \pi)$.
Denote by
\[
{\cal I}: \, \, \fk \, \stackrel{\sim}{\lrw} \, (\fa + \fn)^*
\]
the identification via the pairing $2{\rm Im} \ll \, , \, \gg$, 
and use the same letter to denote the induced map
\[
{\cal I}: \, \, \chi^q(K)^K 
\, \cong \, \wedge^q \fk \,\stackrel{\sim}{\lrw}\, \wedge^q(\fa + \fn)^*,
\]
where $\chi^q(K)^K \cong \wedge^q \fk$ is the 
space of left invariant $q$-vector
fields on $K$.
Let
\[
d_{{\frak a} + {\frak n}}: \, \, \wedge^q(\fa + \fn)^* \, \lrw \,
\wedge^{q+1}(\fa + \fn)^*
\]
be the Chevalley-Eilenberg coboundary operator for the Lie algebra
$\fa + \fn$. The following lemma says that
the (left) $K$-invariant Poisson cohomology of $(K, \pi)$
is isomorphic to the Lie algebra cohomology of the
Lie algebra $\fa + \fn$.

\begin{lem}
\label{lem_K}
We have 
\[
\delta_{\pi} \, = \, {\cal I}^{-1} d_{{\frak a} + {\frak n}} {\cal I}
\]
as degree $1$ operators on $\chi(K)^K \cong \wedge \fk$.
\end{lem}

\noindent
{\bf Proof.} This is again a general fact about Poisson Lie groups specialized
to the Poisson  Lie group $(K, \pi)$. To give a direct proof, we first
notice that both operators are derivations of degree $1$ on $\wedge \fk$.
Thus it is enough to show that they are the same on
$\fk$. This can be checked directly.
\qed

\noindent
{\bf Proof of Theorem \ref{thm_ponf-n}.} 
Denote by $\chi(K)^{K;T}$ the space of multi-vector fields
on $K$ that are both left $K$-invariant and right $T$-invariant.
It is invariant under the operator $\delta_{\pi}$. By Lemma 
\ref{lem_K}, we know that 
\begin{equation}
\label{eq_II}
{\cal I}: \, \, (\chi(K)^{K;T}, \, \delta_{\pi}) \, \lrw \, 
\left( (\wedge(\fa + \fn)^*)^T, \, d_{{\frak a} + {\frak n}} 
\right)
\end{equation}
is an isomorphism of cochain complexes.
Now the projection from $K$ to $K/T$ induces
a surjective cochain complex morphism
\[
p_1: \, \, (\chi(K)^{K;T}, \, \delta_{\pi}) \, \lrw \,  
(\chi(K/T)^K, \, \delta_{\pi_{\infty}}).
\]
Similarly, the map
\[
p_2: \, \, \left( (\wedge(\fa + \fn)^*)^T, \, d_{{\frak a} + {\frak n}}  
\right) \, \lrw \, \left( (\wedge \fn^*)^T, \, d_{\frak n} \right)
\]
induced by the inclusion $\fn \hookrightarrow \fa + \fn$ is
a surjective cochain complex morphism. Since
$\iis p_1 = p_2 {\cal I}$, it follows from 
(\ref{eq_II}) that 
\[
\iis: \, \, (\chi(K/T)^K, \, \delta_{\pi_{\infty}}) \, \lrw \, 
\left( (\wedge \fn^*)^T, \, d_{\frak n} \right) 
\] 
is a cochain complex isomorphism. This proves 
(\ref{eq_ponf-d}). 

We now show that (\ref{eq_ponf-partial}) follows from 
(\ref{eq_ponf-d}). Let $\mu_{0}^{'} \, := \, I_{\infty}^{-1} (\mu_0) 
\in \wedge^n \fn$, 
where $n = \dim_{\Bbb R} \fn$. A general fact about
Lie algebra cohomology says that
\begin{equation}
\label{eq_backl-n}
(d_{\frak n} X) \backl \, \mu_{0}^{'} \, = \, (-1)^{q-1}
b_{\frak n} (X \backl \, \mu_{0}^{'})
\end{equation}
for any $X \in \wedge^q \fn^*$. Since $\mu_{0}^{'}$ is $T$-invariant, the map
$\backl \, \mu_{0}^{'}$ from $\wedge^q \fn$ to
$\wedge^{n-q} \fn^*$ is $T$-equivariant. Thus it induces
an isomorphism from 
$(\wedge^q \fn^*)^T$ to $(\wedge^{n-q} \fn)^T$. Now 
(\ref{eq_ponf-partial}) follows from (\ref{eq_ponf-d})
by comparing
(\ref{eq_pinf-ponf}) and (\ref{eq_backl-n}).
\qed

\begin{rem}
\label{rem_simple-1}
{\em
There is a simpler expression for $\ponf$ as an operator
on the space  $\Omega(K/T)^K$ (as opposed to
on all of $\Omega(K/T)$).
Recall the definition
of $\ponf$ in (\ref{eq_pinf-ponf}) and the definition of 
the $r$-matrix $r$ in (\ref{eq_r}). Use 
the same letter $r$ 
to denote the $K$-invariant bi-vector field  on $K/T$
whose value at $e = eT \in K/T$
is $r$ considered as in $\wedge^2(\fk / \ft)$.
Then, if $X \in \chi^q(K/T)$ is $K$-invariant,
we have $[\pinf, X] = -[r, X]$. A proof similar to 
that of Lemma \ref{lem_delta-partial}
shows that $[r, X] \backl \mu_0 
= (-1)^{q-1} (i_{r} d - 
d i_{r}) (X \backl \mu_0)$. Thus we have
\[
\ponf \, = \, -i_{r} d \, + \, d i_{r}
\]
on $\Omega(K/T)^K$.
}
\end{rem}

\bigskip
We now look at the complexification of
the mixed complex $(\Omega(K/T)^K, \, d, \, \ponf)$. We 
first need the complexification of the Manin
triple $(\fg, \fk, \fa + \fn, 2{\rm Im} \ll \, , \, \gg)$.
For this purpose, we use $J_0$ to denote the complex structure
on $\fg$ as the complexification of $\fk$, and let
$c_0: \fg \rightarrow \fg$ be the complex conjugation
on $\fg$ defined by $\fk$. 
Then the map
\[
\Psi: \, (\fg_{\Bbb C}, \, i) \lrw (\fg \oplus \fg, \, J_0 \oplus J_0): \,
x_1 + i x_2 \Map (c_0(x_1 - J_0 x_2), \, \, x_1 + J_0 x_2)
\]
is an isomorphism of complex Lie algebras. Accordingly, the
complexifications of the various (real) Lie subalgebras of
$\fg$ can be identified via $\Psi$ with (complex)
Lie subalgebras  of $\fg \oplus \fg$.
For example, denoting by
\[
\fg_{1}^{\Delta} \, : = \, \{\, (x, \, x): \, \, x \in \fg_1 \}
\, \subset \fg \oplus \fg
\]
for any complex subspace $\fg_1$ of $\fg$,  since for $x_1, x_2 \in \fk$,
\[
\Psi(x_1 + ix_2) \, = \, (x_1 + J_0 x_2, \, \, x_1 + J_0 x_2),
\]
we have the isomorphism
\[
\Psi: \, \fk_{\Bbb C} \, \stackrel{\sim}{\lrw} \,
\fg^{\Delta}.
\]
Similarly, for $H_1, H_2 \in \fa$,
\[
\Psi(H_1 + i H_2) \, = \, (-H_1 -J_0 H_2, \, \,
H_1 + J_0 H_2),
\]
and thus we have the isomorphism
\[
\Psi: \, \fa_{\Bbb C} \, \stackrel{\sim}{\lrw} \,  \{(-H, \, H): \,
H \in \fh \}.
\]
For $\fn$ as a real Lie algebra, we have
the isomorphism
\[
\Psi: \, \fn_{\Bbb C} \, \stackrel{\sim}{\lrw} \,
\fn_{-} \oplus \fn_{+},
\]
where $\fn_{+} = \fn$ and $\fn_{-} = c_0(\fn)$,
with
\[
\Psi(-{1 \over 2}(\ea + i J_0 \ea)) \, = \, (\eb, \, 0) \hspace{.4in}
\Psi({1 \over 2}(\ea - iJ_0 \ea)) \, = \, (0, \, \ea).
\]
Set
\[
\fb_{-} \oplus_{\frak h} \fb_{+}
\, = \, \{(x_{-} -H, \,\, H + x_{+}): \, x_{\pm} \in
\fn_{\pm}, H \in \fh \}.
\]
Then we have
\[
\Psi: \, (\fa + \fn)_{\Bbb C} \stackrel{\sim}{\lrw} \,
\fb_{-} \oplus_{\frak h} \fb_{+}.
\]
Denote by $\la \, , \ra_{\Bbb C}$ the complex-linear extension of
the bilinear form $\la \, , \, \ra = 2{\rm Im}\ll \, , \gg$
from $\fg$ to $\fg_{\Bbb C}$. Then under the identification of
$\fg_{\Bbb C}$ with $\fg \oplus \fg$ by $\Psi$, it becomes
\[
\la (x_1, x_2), \, (y_1, y_2) \ra_{\Bbb C} \,   =  \,
i(\ll x_1, y_1\gg \, - \,
\ll x_2, y_2 \gg).
\]
The two Lie subalgebras $\fg^{\Delta}$ and $\fb_{-} \oplus_{\frak h} \fb_{+}$
are both isotropic with respect to $\la \, , \ra_{\Bbb C}$.
The triple $(\fg \oplus \fg, \,\, \fg^{\Delta}, \,\,
\fb_{-} \oplus_{\frak h} \fb_{+} )$, together with
$\la \, , \ra_{\Bbb C}$, is a Manin triple of complex Lie
algebras. It is isomorphic to the complexification
$(\fg_{\Bbb C}, \, \,  \fk_{\Bbb C}, \, \, (\fa + \fn)_{\Bbb C}, \, \,
\la \, , \, \ra_{\Bbb C})$
of the Manin triple $(\fg, \, \fk, \, \fa + \fn, \,  \la \, , \, \ra)$
via the map $\Psi$.

\bigskip
Using the above complexification, we have the natural
identification
\[
(\te)_{\Bbb C} \, \cong \, (\fk / \ft)_{\Bbb C} \, \cong \, \fg / \fh.
\]
Using $\fn_{\Bbb C} \cong \fn_{-} \oplus \fn_{+}$ given by $\Psi$,
the complexification of the map $\ii$ in (\ref{eq_ii}),
which we will still denote by $\ii$, becomes
\begin{eqnarray}
\label{eq_I0}
I_{\infty}: & &  \fn_{-} \oplus \fn_{+} \lrw 
(T^{*}_{e}(K/T))_{\Bbb C}:\\
\nonumber
& & \left( I_{\infty}(x_{-}, \, x_{+}), \, \, x \right) \, = \, i
\ll x, \, \, x_{-}\,  - \,x_{+} \gg \, = \,
 \ll x, \, \, ix_{-} \, - \, i x_{+} \gg
\end{eqnarray}
for $x \in \fg / \fh \cong (T_{e}(K/T))_{\Bbb C}$. 

Notice the minus sign in front of $ix_{+}$. Notice also
that there is another more standard identification
\begin{equation}
\label{eq_I}
I: \, \,  \fn_{-} \oplus \fn_{+} \lrw (T^{*}_{e}(K/T))_{\Bbb C}:
\, \, 
 \left( I(x_{-}, \, x_{+}), \, \, x \right) \, = \, 
\ll x, \, \, x_{-} \,  + \, x_{+} \gg 
\end{equation}
for $x \in \fg / \fh \cong (T_{e}(K/T))_{\Bbb C}$. As we will see,
it is very important that our identification $\ii$  is different
from the standard one $I$.

\bigskip
At this point, we introduce a complex structure $J$ on $K/T$:
 at $\eo \in K/T$,
it is the complex structure on $\te \cong \fk / \ft$ coming from the
identification $\fk /\ft \cong \fg / \fb_{-}$, where 
$\fb_{-} = \fn_{-} + \fh$. 
We extend $J$ to all of $K/T$ by left translations by elements in $K$.
We will use the same letter $J$ to denote its dual map
on $T^*(K/T)$ as well as its complex linear extension to
the complexified cotangent bundle $T^{*}_{\Bbb C}(K/T)$ and
its multi-linear extension to $\wedge T^{*}_{\Bbb C}(K/T)$.
If we identify
$\fk / \ft \cong {\rm span}_{\Bbb R} \{\Xa, \Ya: \a > 0\}$, then
\[ 
J \Xa \, = \, \Ya, \hspace{.3in} J \Ya \, = \, - \Xa. 
\] 
Thus, if we use $\{X^{\alpha}, Y^{\alpha}: \, \a > 0 \}$ to denote
the dual basis for $\tse$, then 
\[ 
J X^{\alpha} \, = \, - Y^{\alpha}, \hspace{.3in}  
J Y^{\alpha} \, = \, X^{\alpha}. 
\] 

\begin{lem}
\label{lem_I-J}
1) Under the identifications $\ii$ in (\ref{eq_I0})
or $I$ in (\ref{eq_I}), the operator $J$ on 
$(\tse)_{\Bbb C}$ becomes the following one on $\fn_{-} \oplus \fn_{+}$:
\[
I^{-1} J I \, = \, I_{\infty}^{-1} J \ii: \, \, 
\fn_{-} \oplus \fn_{+} \lrw \fn_{-} \oplus \fn_{+}: \, \,
(x_{-}, \, x_{+}) \Map (i x_{-}, \, \, -i x_{+}).
\]

2) The identifications $\ii$ and $I$ are related by
\[
\ii I^{-1} \, = \, J: \, \, 
(\tse)_{\Bbb C} \lrw (\tse)_{\Bbb C}.
\]
\end{lem}

\noindent
{\bf Proof.} 
1) is straightforward, and 2) follows from the definitions
of $\ii$ and $I$ and from 1).
\qed

\bigskip
Denote by $C = \Omega(K/T, \C)^K$ the space of all $K$-invariant
complex valued differential forms on $K/T$, and by $C^q =
\Omega^q(K/T, \C)^K$ those of degree $q$.
Denote the  multilinear extension of $I_{\infty}$ to
$(\wedge^{\bullet}(\fn_{-} \oplus \fn_{+}))^T$ by the same letter:
\begin{equation}
\label{eq_ii-multi}
I_{\infty}: \, \, (\wedge^{q}(\fn_{-} \oplus \fn_{+}))^T
\, \stackrel{\sim}{\lrw} \, C^{q}.
\end{equation}
By Theorem \ref{thm_ponf-n}, we have

\begin{thm}
\label{thm_ponf-n-complex}
As operators on $C$, 
\begin{equation}
\label{eq_ponf-n-complex}
\ponf \, = \, \ii b_{{\frak n}_{-} \oplus {\frak n}_{+}}
I_{\infty}^{-1}: \, \, C^q \lrw C^{q-1},
\end{equation}
where $ b_{{\frak n}_{-} \oplus {\frak n}_{+}}$ is the Chevalley-Eilenberg
 boundary operator for the Lie algebra $\fn_{-} \oplus \fn_{+}$
as well as it restriction to
$\wedge(\fn_{-} \oplus \fn_{+})^T$.
\end{thm}

\bigskip
The complex structure
$J$ on $K/T$ defines a bi-grading on $C$:
\[
C \, = \, \oplus_{p,q} C^{p,q}
\]
where $C^{p,q}$ is the space of $K$-invariant differential forms on $K/T$ of
holomorphic degree $p$ and anti-holomorphic degree $q$.
By Lemma \ref{lem_I-J}, 
\[
\ii  \, \, {\rm and} \, \, I: \, \, (\wedge^p \fn_{-} \otimes 
\wedge^q \fn_{+})^T \,
 \lrw \, C^{p,q}
\]
are both isomorphisms.

\bigskip
For the operators $d$ and $\ponf$, denote by the
same letters their complex linear extensions to $C$.
Since $J$ is integrable, we can write $d$ as
\[
d \, = \, d' \, + \, d^{''},
\]
where $d^{'}:  C^{p,q} \lrw C^{p+1, q}$ and $ d^{''}: C^{p,q} 
\lrw C^{p, q+1}$.

\bigskip
We now look at how the operator $\ponf$ respects the  
 bi-grading of $C$.
Denote by
$\Pi_{p.q}$ the projection from $(\wedge \fn_{-} \otimes 
\wedge \fn_{+})^T$ to $(\wedge^p \fn_{-} \otimes 
\wedge^q \fn_{+})^T$, and set
\[
\tau \, = \, \sum_{p,q} (-1)^p \Pi_{p,q}.
\]
Since $\fn_{-} \oplus \fn_{+}$ has the direct sum Lie algebra
structure, we have
\[
b_{{\frak n}_{-} \oplus {\frak n}_{+}} \, = \, 
b_{{\frak n}_{-}} \otimes 1 \, + \, \tau (1 \ot b_{{\frak n}_{+}}),
\]
where $b_{{\frak n}_{-}}$ and $b_{{\frak n}_{+}}$ are 
the Chevalley-Eilenberg boundary operators for the Lie
algebras $\fn_{-}$ and $\fn_{+}$ respectively. Set
\beqa
& & b^{'} \, = \, b_{{\frak n}_{-}} \otimes 1: \, \, (\wedge^p \fn_{-} 
\otimes  
\wedge^q \fn_{+})^T \lrw (\wedge^{p-1} \fn_{-} \otimes  
\wedge^q \fn_{+})^T \\
& & b^{''} \, = \, \tau (1 \ot b_{{\frak n}_{+}}): \, \,
(\wedge^p \fn_{-} \otimes  
\wedge^q \fn_{+})^T \lrw (\wedge^p \fn_{-} \otimes  
\wedge^{q-1} \fn_{+})^T\\
& &\partial_{\infty}^{'} \, = \, \ii b^{'} I_{\infty}^{-1}: \, \, C^{p,q} 
\lrw C^{p-1,q} \\
& & \partial_{\infty}^{''} \, = \, \ii b^{''} I_{\infty}^{-1}: \, \, 
C^{p,q} \lrw C^{p, q-1}.
\eeqa
Then, by Theorem \ref{thm_ponf-n-complex}, we have
\[
\ponf \, = \, \partial_{\infty}^{'} \, + \, \partial_{\infty}^{''}.
\]
Recall that the two operators $d$ and $\ponf$ anti-commute on $C$, i.e.,
\[
d \ponf \, + \, \ponf d \, = \, 0.
\]
Thus by reasons of degree, we have

\begin{prop}
\label{prop_primes}
The following holds on $C = \Omega(K/T, \C)^K$:
\beqa
& & (d^{'} \partial_{\infty}^{'} \, + \, 
\partial_{\infty}^{'} d^{'} ) \, + \, (d^{''} \partial_{\infty}^{''}
\, + \, \partial_{\infty}^{''} d^{''}) \, = \, 0;\\
& & d^{'} \partial_{\infty}^{''} \, + \, \partial_{\infty}^{''} 
d^{'} \, = \, 0;\\
& & d^{''} \partial_{\infty}^{'} \, + \, \partial_{\infty}^{'}  
d^{''} \, = \, 0.
\eeqa
\end{prop}

\begin{rem}
\label{rem_can-also-be-proved}
{\em
By Remark \ref{rem_simple-1}, we also have
\beqa
\partial_{\infty}^{'} & = & -i_r d^{''} \, + \, d^{''} i_r \\
\partial_{\infty}^{''} & = & -i_r d^{'} \, + \, d^{'} i_r.
\eeqa
The statements in Proposition \ref{prop_primes} also follow from these 
two identities.
}
\end{rem}

\subsection{Kostant's operator $\partial$ versus the Koszul-Brylinski
operator $\ponf$}
\label{sec_partial}

We now recall the operator $\partial$ on $C = \Omega(K/T, \C)^K$ 
introduced 
by Kostant in \cite{ko:63}.

\bigskip
In \cite{ko:63}, the map $I$ given in (\ref{eq_I}) was used to identify
$\fn_{-} \oplus \fn_{+}$ with $(\tse)_{\Bbb C}$. Set
\begin{equation}
\label{eq_partial-ko}
\partial \, = \, - \,I b_{{\frak n}_{-} \oplus
{\frak n}_{+}} I^{-1}: \, \, C^q \lrw C^{q-1}.
\end{equation}
This is the operator $\partial$
in \cite{ko:63}. We will call $\partial$ the {\bf Kostant operator}.
The relation between the Kostant operator $\partial$ and the 
Koszul-Brylinski operator $\ponf$ defined by the Bruhat Poisson 
structure $\pinf$ is now clear from 
(\ref{eq_ponf-n-complex}) and Lemma \ref{lem_I-J}:

\bigskip
\begin{thm}
\label{thm_main-1}
The Kostant operator $\partial$ and the Koszul-Brylinski operator
$\ponf$ on $C = \Omega(K/T, \C)^K$ are related by
\[
\partial \, = \, - J^{-1} \ponf J \, = \, J \ponf J^{-1},
\]
where $J$ is the complex structure on $K/T$ (defined before Lemma
\ref{lem_I-J}).
\end{thm}

\bigskip
\begin{cor}
\label{cor_1}
On $C = \Omega(K/T, \C)^K = \oplus_{p,q} C^{p,q}$, we have
\[
\partial \, = \, - i\partial_{\infty}^{'} \, + \, i \partial_{\infty}^{''}.
\]
\end{cor}

\noindent
{\bf Proof.} This follows immediately from Theorem \ref{thm_main-1}
and the fact that $J|_{C^{p,q}} = i^{p-q} {\rm id}$.
\qed

\bigskip
Set
\[
\partial^{'} \,  = \, - \,i \partial_{\infty}^{'}, \hspace{.4in}
\partial^{''} \,  = \,  i \partial_{\infty}^{''},
\]
so that $ \partial  =  \partial^{'}  +  \partial^{''}$.
By Proposition \ref{prop_primes} and Corollary \ref{cor_1}, we have

\begin{cor}
\label{cor_2}
The following hold on $C = \Omega(K/T, \C)^K = \oplus_{p,q} C^{p,q}$:
\beqa
& & d^{'} \partial^{'} \, + \, \partial^{'} d^{'} \, = \, 
 d^{''} \partial^{''} \, + \, \partial^{''} d^{''};\\
& &  d^{'} \partial^{''} \, + \, \partial^{''} d^{'} \, = \, 0;\\
& &  d^{''} \partial^{'} \, + \, \partial^{'} d^{''} \, = \, 0.
\eeqa
\end{cor}

\begin{rem}
\label{rem_kos}
{\em The above facts are in Proposition 4.2 of 
\cite{ko:63}.
}
\end{rem}

\bigskip
\begin{thm}
\label{thm_main-2}
A form  $\xi \in C^{p,q}$ is $\partial$-closed if and only if it is 
$\ponf$-closed.
\end{thm}
 
\noindent
{\bf Proof.} Let $\xi \in C^{p,q}$. Then
\[
\partial \xi \, = \, -\, i  \partial_{\infty}^{'} \xi \, + \, 
i \partial_{\infty}^{''} \xi
\hspace{.3in} {\rm and} \hspace{.3in}
\partial_{\infty} \xi \, = \, \partial_{\infty}^{'} \xi \, + \,
\partial_{\infty}^{''} \xi.
\]
Thus $\partial \xi = 0$ if and only if $\partial_{\infty}^{'} \xi = 0$ and
$\partial_{\infty}^{''} \xi = 0$, or, if and only if $\ponf \xi = 0$.
\qed

\bigskip
In \cite{ko:63}, a $K$-invariant form $\xi$ on $K/T$ is said
to be harmonic if it satisfies $d \xi = \partial \xi = 0$. 
We will also say that such forms are $(d, \partial)$-harmonic
or say that they are ``harmonic in the sense of Kostant".

\bigskip
\begin{cor}
\label{cor_main-2}
A $K$-invariant form $\xi$ on $K/T$ of pure bi-degree is
harmonic in the sense of Kostant if and only if it is harmonic
 with respect to the Bruhat Poisson structure $\pinf$ (and a
$K$-invariant volume form $\mu_0$) on $K/T$.
\end{cor}

\subsection{Kostant's Theorems}
\label{sec_kos-thm}

In this section, we will recall the main theorems in
\cite{ko:63}. We will give new proofs of some of them using the Bruhat
Poisson structure. 

\bigskip
\noindent
{\bf Kostant's Theorem 1}  (Theorem 4.5 in \cite{ko:63}) 
{\it The operators $d$ and $\partial$ on $C = \Omega(K/T, \C)^K$
are disjoint, i.e., 
\[
{\rm Im} (d )\cap {\rm Ker} (\partial) \, = \, {\rm Im} 
(\partial) \cap
{\rm Ker} (d) \, = \, 0.
\]
}

\bigskip
In Section \ref{sec_another-proof}, we will give a proof of this theorem 
by using the Koszul-Brylinski operators for a family of symplectic 
structures on $K/T$. In this section,
we assume this theorem and proceed
to prove the other main results in \cite{ko:63}
using our Poisson
interpretation of the operator $\partial$.

\bigskip
Set
\[
S \, = \, d \partial \, + \, \partial d \, \in \, {\rm End}(C).
\]
It follows immediately from the disjointness of
$d$ and $\partial$ that $\xi \in {\rm Ker} (S)$ if and only if
$d \xi = \partial \xi = 0$, i.e., $\xi$ is $(d, \partial)$-harmonic,
and that the maps
\begin{eqnarray}
\label{eq_fd}
\psi_{d, S}: \, \, {\rm Ker} (S) \lrw H(C, d): \, \, \xi \Map [\xi]_d,
\\
\label{eq_fp}
\psi_{\partial, S}: \, \, {\rm Ker} (S) \lrw H(C, \partial): \, \, \xi 
\Map [\xi]_{\partial}
\end{eqnarray}
are isomorphisms, where $H(C, d)$ and $H(C, \partial)$ are the
cohomology groups of $C$ defined by $d$ and $\partial$,
and $[\xi]_d$ and $[\xi]_{\partial}$
are the cohomology classes defined by $\xi$ in  
$H(C, d)$ and $H(C, \partial)$ respectively. 

It follows from
Corollary \ref{cor_2} that
\[
S \, = \, 2(d^{'} \partial^{'} \, + \, 
\partial^{'} d^{'}) \, = \, 
2(d^{''} \partial^{''} \, + \, 
\partial^{''} d^{''}): \, \, C^{p,q} \lrw C^{p,q}
\]
(this is Proposition 4.2 in \cite{ko:63}).
Consequently, all the three spaces 
${\rm Ker} (S), \, 
H(C, d)$ and $H(C, \partial)$ are
bi-graded,  and the maps $\psi_{d, S}$ and $\psi_{\partial, S}$ are 
isomorphisms of bi-degree $(0,0)$. 
 The
space $H(C, d)$ is of course the de Rham cohomology of
$K/T$ (with complex coefficients), and it is clear from the
definition of $\partial$ that
\[
I: \, \, (H_{\bullet}(\fn_{-}) \ot H_{\bullet}(\fn_{+}))^T 
\lrw H(C, \partial)
\]
is a bi-degree $(0,0)$ isomorphism, where
$H_{\bullet}(\fn_{-})$ and $H_{\bullet}(\fn_{+})$
are respectively the Lie algebra homology groups of
the Lie algebras $\fn_{-}$ and $\fn_{+}$.

\bigskip
The space $(H_{\bullet}(\fn_{-}) \ot H_{\bullet}(\fn_{+}))^T$
has a distinguished basis.
Indeed, the space 
$H_{\bullet}(\fn_{+})$ as a $T$-module is shown in 
\cite{ko:61} to have weights exactly of the form
$\rho - w \rho$, where $w$ is in the Weyl group $W$,
and each such weight has multiplicity $1$.  
By pairing $\fn_{-}$ and $\fn_{+}$ by the Killing form $\ll \, , 
\, \gg$, we can regard $H_{\bullet}(\fn_{-})$
as the contragradient of $H_{\bullet}(\fn_{+})$
as $T$-modules. The weight decomposition of
$H_{\bullet}(\fn_{+})$ then gives a
basis for
$(H_{\bullet}(\fn_{-}) \ot H_{\bullet}(\fn_{+}))^T$ by Schur's Lemma.
Explicitly, for $w \in W$, let 
\[
\Phi_w \, = \, \{\alpha_1, \alpha_2, ..., \alpha_{l(w)} \} \, = \, 
\{\alpha > 0: \, w^{-1} \alpha < 0 \},
\]
where $l(w)$ is the length of $w$. Set
$\beta_j =  - w^{-1} \a_j$ for $j = 1, .., l(w)$.
Then
\[
\{\beta_1, \, \beta_2, \, ..., \,\beta_{l(w)} \} \, = \, \Phi_{w^{-1}}.
\]
Set
\begin{equation}
\label{eq_h}
h^{w^{-1}} ~ = ~ (-1)^{{\frac{l(w)(l(w)-1)}{2}}} 
i^{l(w)} ~ E_{-\beta_1} \wedge E_{-\beta_2}
\wedge \cdots \wedge E_{-\beta_{l(w)}} ~ \ot ~
E_{\beta_1} \wedge E_{\beta_2} \wedge \cdots \wedge 
E_{\beta_{l(w)}}.
\end{equation}
Then $b_{{\frak n}_{-} \oplus {\frak n}_{+}} (h^{w^{-1}}) = 0$,
and, up to scalar multiples, 
the set $\{[h^{w^{-1}}]: w \in W \}$ is the basis of 
$(H_{\bullet}(\fn_{-}) \ot H_{\bullet}(\fn_{+}))^T$
determined by the weight decomposition of $H_{\bullet}(\fn_{+})$.

\bigskip
Consider now the $K$-invariant forms $I(h^{w^{-1}})
\in C^{l(w), l(w)}$. They are $\partial$-closed and
thus also $\ponf$-closed by Theorem
\ref{thm_main-2}. The following lemma will be used later.

\bigskip
\begin{lem}
\label{lem_h}
Let $w, w_1 \in W$ be such that
$l(w) = l(w_1)$. Regard $w$ and $w_1$ as points in $K/T$.
Then
\begin{equation}
\label{eq_h-p}
\left( I(h^{w_{1}^{-1}}), \, \, \, {\pi_{\infty}^{l(w_1)}
\over l(w_1)!} \right) (w) \, = \, \delta_{w, w_1},
\end{equation}
where $\pi_{\infty}^{l(w_1)}$ is the wedge product of 
$\pi_{\infty}$ with itself $l(w_1)$-times.
\end{lem}

\noindent
{\bf Proof.} Since $I(h^{w_{1}^{-1}})$ is $K$-invariant, we 
need to find the value 
\[
L_{{\dot{w}}^{-1}} \pi_{\infty}^{l(w)} (w) \, \in \, \wedge^{2l(w)} \te,
\]
where $\dot{w} \in K$ is any representative
of $w$ in $K$, and $L_{\dot{w}^{-1}}$ is
the map $K/T \rightarrow K/T: kT \mapsto \dot{w}^{-1} kT$ as well
as its differential. By the definition of the Poisson structure
$\pi$ on $K$, we easily see that
\[
L_{{\dot{w}}^{-1}} \pi(\dot{w}) \, = \, 
-i (E_{\beta_1} \wedge E_{-\beta_1} \, + \, \cdots \, + \, 
E_{\beta_{l(w)}} \wedge E_{-\beta_{l(w)}}) \in \fg \wedge \fg.
\]
Here $L_{{\dot{w}}^{-1}}$ also denotes the left translation
on $K$ by $\dot{w}^{-1}$.
Hence
\[
{L_{{\dot{w}}^{-1}} \pi_{\infty}^{l(w)} (w) \over l(w)!} \, = \, 
(-i)^{l(w)} (-1)^{l(w) (l(w) -1) \over 2} 
E_{-\beta_1} \wedge E_{-\beta_2}
\wedge \cdots \wedge E_{-\beta_{l(w)}} ~ \ot ~
E_{\beta_1} \wedge E_{\beta_2} \wedge \cdots \wedge
E_{\beta_{l(w)}}.
\]
Now (\ref{eq_h-p}) follows from the definitions of $I$ and 
of $h^{w_{1}^{-1}}$.
\qed

The differential forms $I(h^{w^{-1}})$, for
$w \in W$, although $\ponf$-closed,
are in general not $d$-closed, and thus not harmonic
(in the sense of Kostant or with respect to $\pinf$). But recall
that the map $\psi_{\partial, S}: {\rm Ker} (S) \rightarrow H(C, \partial)$
in (\ref{eq_fp}) is a bi-degree $(0, 0)$ isomorphism. It is used 
in \cite{ko:63} to define harmonic forms:

\bigskip
\begin{dfn}
\label{dfn_kostant-forms}
{\em
For $w \in W$, set
\[
s^w \, = \, \psi_{\partial, S}^{-1} \left([I(h^{w^{-1}})]_{\partial}
\right) \, \in \, {\rm Ker} (S).
\]
These are called the {\bf Kostant harmonic forms}.
}
\end{dfn}

\bigskip
Thus, by definition, the form $s^w$ has bi-degree $(l(w), l(w))$,
and it is $(d, \partial)$ as well as $(d, \ponf)$-harmonic.
In other words, the Kostant harmonic forms are Poisson harmonic.

\bigskip
We now review the explicit formulas for $s^w, w \in W$, as given in 
\cite{ko:63}. 

\bigskip
Introduce the operators $E_I$ and $L_0$ on $(\wedge \fn_{-}
\ot \wedge \fn_{+})^T$ by
\begin{equation}
\label{eq_E}
E_I ~ = ~ 2 \sum_{\alpha > 0} ad_{E_{-\alpha}} ~ \ot ~
ad_{E_{\alpha}}
\end{equation}
and
\beqa
& ~~~ & L_0|_{{\Bbb C} E_{-\alpha_1} \wedge E_{-\alpha_2}
\wedge \cdots \wedge E_{-\alpha_{p}}
\ot  E_{\beta_1} \wedge E_{\beta_2}  \wedge \cdots \wedge
E_{\beta_{q}}} \vspace{.2in} \\
 \hspace{.2in}  &  = & \left\{
\begin{array}{ll} 0 & {\rm if} ~ \| \rho \|^2 - \|\rho -(\alpha_1
+ \alpha_2 + \cdots + \alpha_p) \|^2 = 0 \vspace{.1in} \\
{\frac{1}{ \| \rho \|^2 - \|\rho -(\alpha_1
+ \alpha_2 + \cdots + \alpha_p) \|^2}} {\rm id}& {\rm if} ~
 \| \rho \|^2 - \|\rho -(\alpha_1
+ \alpha_2 +  \cdots + \alpha_p) \|^2 \neq 0 \end{array} \right.,
\eeqa
where $\| \lambda \|^2 = \ll \lambda, \lambda \gg$ for $\lambda \in \fh^*$.
Set \[
R ~ = ~ - L_0 E_I.
\]

\bigskip
\noindent
{\bf Kostant's Theorem 2} (Theorem 5.6 in \cite{ko:63})
The harmonic forms $s^w$, for $w \in W$,  are explicitly given by
\begin{equation}
\label{eq_s}
s^w ~ = ~ I \left( (1 - R)^{-1} h^{w^{-1}} \right) ~ = ~
I(h^{w^{-1}} ~ + ~
R h^{w^{-1}} ~ + ~ R^2 h^{w^{-1}} ~ + ~ \cdots),
\end{equation}
where $h^{w^{-1}}$ is given by (\ref{eq_h}).

\bigskip
The proof of this theorem is relatively easy, and it uses a fact
that we will need later, so we review the main points in the proof now.

\bigskip
We introduce a Hermitian inner product $h$ on $C = \Omega(K/T, \C)^K$
as follows: Identify 
\[
\te \, \cong \, \fk / \ft \, \cong \, {\rm span}_{\Bbb R}
\{ \Xa, \, \Ya: \, \a > 0 \}
\]
and consider the symmetric and positive definite scalar product 
$g^{-1}$ on $\te$ given by
\[
g^{-1}(\Xa, \, \Xa) \, = \, g^{-1}(\Ya, \, \Ya) \, = \,  2,
\hspace{.3in} \forall \a > 0
\]
and $0$ otherwise. Denote by $\{X^{\alpha},  Y^{\alpha}:  \a > 0 \}$
the basis of $\tse$ that is dual to the basis of 
$\te$ given by $\{X_{\alpha},  Y_{\alpha}: 
\a > 0 \}$. Then $g^{-1}$ induces the following scalar product $g$ on
$\tse$:
\begin{equation}
\label{eq_g}
g(X^{\alpha}, \, X^{\alpha}) \, = \, g(Y^{\alpha}, \, Y^{\alpha}) 
\, = \, {1 \over 2}, \hspace{.3in} \forall \a > 0 
\end{equation} 
and $0$ otherwise. Now define $h$ to be the Hermitian extension of $g$
from $\tse$ to $(\tse)_{\Bbb C}$. It is easy to see that under the 
identification $I$ given in (\ref{eq_I}), the induced Hermitian
product on $\fn_{-} \oplus \fn_{+}$ is given by
\begin{equation}
\label{eq_h-I}
(h \circ I) (\ea, \, \ea) \, = \, (h \circ I)(\eb, \, \eb) \, = \, 1
\hspace{.3in} \forall \a > 0
\end{equation}
and $0$ otherwise. Notice that, since the complex structure $J$
is an isometry for $h$, we have $h \circ I = h \circ \ii$ by Lemma
\ref{lem_I-J}.

Let $\delta$ be the adjoint operator of $\partial$ with respect to
the Hermitian product $h$. Set
\[
L \, = \, \partial \delta \, + \, \delta \partial: \, \, C^q \lrw C^q.
\]
The eigenvalues and eigenvectors of this Laplacian $L$
were determined by Kostant in \cite{ko:61}.  They are easier to 
describe if we identify $C$ with $\wedge(\fn_{-} \oplus \fn_{+})^T$:
For notational simplicity, set $L_{I} = I^{-1} L I$, so it
is an operator on $\wedge(\fn_{-} \oplus \fn_{+})^T$. Kostant showed in 
\cite{ko:61} that with respect to the basis of 
$\wedge(\fn_{-} \oplus \fn_{+})^T$ given by
\[
\{ E_{-(\alpha)} \otimes E_{(\beta)}: \, \, 
(\a) = \{ \a_1,..., \a_p \}, \, (\beta) = \{\beta_1, ..., \beta_q \}:
\, \alpha_1 + \cdots + \alpha_p = \beta_1 + \cdots + \beta_q\},
\]
where
\[
E_{-(\alpha)} \, = \, E_{-\alpha_1} \wedge \cdots \wedge E_{-\alpha_p}
\hspace{.2in} {\rm and} \hspace{.2in}
E_{(\beta)} \, = \, E_{\beta_1} \wedge \cdots \wedge E_{\beta_q},
\]
the Laplacian $L_I$ is diagonal, and that
$L_I(E_{-(\alpha)} \otimes E_{(\beta)}) = 0$ if and only if 
$(\alpha) = \Phi_w = \{ \a > 0: w^{-1} \a < 0 \}$
for some $w \in W$ and otherwise
\[
L_I(E_{-(\alpha)} \otimes E_{(\beta)}) \, = \, 
(\| \rho \|^2 - \|\rho -(\alpha_1
+ \alpha_2 + \cdots + \alpha_p) \|^2) 
E_{-(\alpha)} \otimes E_{(\beta)}.
\]
Thus the operator $L_0$ is nothing but the Green's operator
for $L_I$.

\bigskip
Now let $E \in \End(C)$ be such that
\begin{equation}
\label{eq_S-L}
S \, = \, L \, + \, E.
\end{equation}
Let $E_I = I^{-1} E I$ be the corresponding operator on
$\wedge(\fn_{-} \oplus \fn_{+})^T$.

\begin{lem}
\label{lem_E}
The operator $E_I$ is explicitly given by (\ref{eq_E}).
\end{lem}

This is Proposition 2.8 in \cite{ko:63} (applied to
the case when ${\frak r} = \fn_{-} + \fn_{+}$ in the notation of
\cite{ko:63}). The proof as given in \cite{ko:63} is very simple,
and we omit it here.

\bigskip
The explicit formula for the form $s^w$ given in
(\ref{eq_s}) is now very easy to derive: By definition,
the element $h^{w^{-1}}$ is the projection of $I^{-1}(s^w)$
to $\Ker (L_I)$ with respect to the decomposition
$\wedge(\fn_{-} \oplus \fn_{+})^T = \Im (L_I) + \Ker (L_I)$.
Thus
\beqa
I^{-1}(s^w) &  =  & h^{w^{-1}} \, + \, L_I L_0 (I^{-1}(s^w)) \\
& = & h^{w^{-1}} \, + \, L_0 L_I (I^{-1}(s^w)) \\
& = & h^{w^{-1}} \, - \, L_0 E_I (I^{-1}(s^w)).
\eeqa
Hence $s^w = I(1 + L_0 E_I)^{-1} (h^{w^{-1}})$.

\begin{rem}
\label{rem_p-commute-E-L}
{\em
It is clear that $\partial$ commutes with $L$ and thus also with
the Green's operator of $L$ on $C$. Since $\partial$ 
commutes with $S$, it also commutes with $E = S - L$.
This fact will be used in 
Section \ref{sec_another-proof}.
}
\end{rem}

\begin{exam}
\label{exam_not-adjoint}
{\em
The operator $S$ on $C$ is in general not semi-simple, as is shown in
the following example. Consider the case when $\fg = sl(4, \C)$.
Consider the space spanned by the three vectors
\[
x \, = \, E_{23} \wedge E_{14} \ot E_{32} \wedge E_{41}
\hspace{.2in} 
y \, = \, E_{13} \wedge E_{14} \ot E_{31} \wedge E_{41}
\hspace{.2in}
z \, = \, E_{24} \wedge E_{14} \ot E_{42} \wedge E_{41}
\]
in $(\wedge^2 \fn_{-} \ot \wedge^2 \fn_{+})^T$. Denote by  $S_I$ 
the operator
$I^{-1} S I$ corresponding to $S.$
 It is easy to check that
the Laplacian $L_I$ acts by the same scalar on $x,$ $y,$ and $z$ while
$E_I$ is strictly upper triangular with respect to the ordered basis
$(x,y,z).$ It follows that $S_I=L_I+E_I$ is not semisimple on the space
spanned by these three vectors, and hence that $S$ is not semisimple
on $C$.  This example shows that the Kostant operator $\partial$ is
in general not adjoint to the de Rham operator $d$
with respect to any Hermitian metric on $K/T$,
for the operator $S$ would be semi-simple if this were the case.
}
\end{exam}

\bigskip
We now state the second main result in \cite{ko:63}, which is
the first two properties of the forms $s^w,
w \in W$ listed in the Introduction.

\bigskip
\noindent
{\bf Kostant's Theorem 3} 
(Theorem 6.15 and Corollary 6.15 in \cite{ko:63}) 
{\it The
de Rham cohomology classes $[s^w]_d$, for $w \in W$, form
a basis of $H(K/T, \C)$ that is, up to scalar multiples $\lambda_w$,
 dual to the 
basis of the homology of $K/T$ formed by the Schubert varieties. 
In fact, when $l(w) = l(w_1)$ but $w \neq w_1$,
the form $i_{w_1}^{*} s^w$ on the Schubert cell $\Sigma_{w_1}$
is identically zero, where $i_{w_1}: \Sigma_{w_1} \hookrightarrow
K/T$ is the inclusion map.
}

\bigskip
The numbers $\lambda_w$ were later calculated by
Kostant and Kumar in \cite{k-k:integral} to be
\begin{equation}
\label{eq_lambda-w}
\lambda_w \, = \, \prod_{j=1}^{l(w)} 
{2 \pi \over \ll \rho, \, \a_j \gg}.
\end{equation}
(Here, unlike elsewhere in the paper, $\pi$ 
denotes the irrational number 3.14159....)

\bigskip
We will now give a proof of the above properties 
using the 
Poisson interpretation of the harmonic forms
$s^w$.

\bigskip
\noindent
{\bf A Poisson Geometrical Proof of Kostant's Theorem 3.}
Since the maps $\psi_{d, S}$ and $\psi_{\partial, S}$ in (\ref{eq_fd})
and (\ref{eq_fp}) are both isomorphisms, we know that
$\{[s^w]_d: \, w \in W \}$ is
a basis for the de Rham cohomology space
$H(K/T, \C)$. It remains to show that
\[
\int_{\Sigma_{w_1}} s^w \, = \, 
\left\{ \begin{array}{ll} 0 & {\rm if} \, \, w_1 \neq w \\
\lambda_w & {\rm if} \, \, w_1 = w \end{array} \right.
\]
where $\lambda_w$ is given in (\ref{eq_lambda-w}).
This is clearly true if $l(w) \neq l(w_1)$.
Assume now that $l(w) = l(w_1)$. For notational 
simplicity, we set
\[
s^{w}_{w_1} \, = \, i_{w_{1}}^{*} s^w \, \in \, 
\Omega^{2l(w_1)} (\Sigma_{w_1}).
\]
Denote by $\Omega_{w_1}$ the symplectic $2$-form
on $\Sigma_{w_1}$ induced from $\pinf$ and by 
$\mu_{w_1} = \Omega_{w_1}^{l(w_1)} / l(w_1)!$
the Liouville volume form 
on $\Sigma_{w_1}$ defined by $\Omega_{w_1}$. 
We now relate $s^{w}_{w_1}$ and $\mu_{w_1}$.

Recall that the Bruhat Poisson structure $\pinf$ is invariant under 
the $T$-action on $K/T$ by left translations. All the 
Schubert cells are invariant under this $T$ action. Thus,
since $\Sigma_{w_1}$ is simply-connected, there is a unique 
moment map, denoted by $\phi_{w_1}: 
\Sigma_{w_1} \rightarrow \ft^*$, 
for the $T$-action on $\Sigma_{w_1}$ such that $\phi_{w_1}(w_1) = 0$.
The function $\la \phi_{w_1}, \, iH_{\rho} \ra$ is then 
a Hamiltonian function for the vector field $\theta_0$ on
the symplectic manifold $(\Sigma_{w_1}, \, \Omega_{w_1})$.
Consider the volume form $\mu = e^{\la \phi_{w_1}, \, iH_{\rho} \ra}
\mu_{w_1}$ on $\Sigma_{w_1}$. Then by Example \ref{exam_2} in Section
\ref{sec_poi-harm}, a top degree form $s$ on $\Sigma_{w_1}$
is harmonic with respect to $\pi|_{\Sigma_{w_1}}$ and $\mu$
if and only if $s = c \mu$ for a constant function $c$ on $\Sigma_{w_1}$.
But $\ponf s^w = 0$ implies that $s^{w}_{w_1}$ is such a form
on $\Sigma_{w_1}$. Thus there exists a constant $C^{w}_{w_1}$
such that
\[
\sww \, = \, \lw e^{\la \phi_{w_1}, \, iH_{\rho} \ra}
\mu_{w_1}.
\]
By Lemma \ref{lem_h} and by the explicit formula
for $s^w$, we know that 
\[
\left(s^{w_1}, \, \, \, {\pi_{\infty}^{l(w_1)} \over
l(w_1)!} \right) (w) \, = \, \delta_{w, w_1}.
\]
Thus $\lw 
= \delta_{w, w_1}$ and 
\[
i_{w_1}^{*} s^w \, = \, 
\left\{ 
\begin{array}{ll} 0 & {\rm if} \, \, 
l(w) = l(w_1), w \neq w_1 \\
e^{\la \phi_{w_1}, \, iH_{\rho} \ra}
\mu_{w_1} & {\rm if} \, \, w = w_1. 
\end{array} \right.
\]
The fact that $\int_{\Sigma_w} s^w  =  \lambda_w$ as
given in (\ref{eq_lambda-w}) now follows from 
the Duistermaat-Heckman formula. See Remark 
\ref{rem_integral}.
\qed

\begin{rem}
\label{rem_integral}
{\em
Explicit formulas for the quantities $s^w, \, 
\mu_w$ and $\la \phi_{w_1}, iH_{\rho} \ra$ are written down in
\cite{lu:coor} in certain coordinates on the Schubert cell $\Sigma_w$.
In particular, the function $\la \phi_{w_1}, iH_{\rho} \ra$
goes to $-\infty$ towards the boundary of the closure of
$\Sigma_w$ in $K/T$. 
The integral $\int_{\Sigma_w} s^w$ can also be calculated 
directly using the formula for $s^w$ in these coordinates.
See \cite{lu:coor}.
}
\end{rem}

In \cite{ko:63}, Kostant actually
constructed harmonic forms on $G/P$ for any parabolic
subgroup $P$ of $G$ that contains $B_{+}$. He showed that
with respect to the projection $ \nu: G/B_{+} \rightarrow G/P$ induced
by the inclusion $B_{+} \rightarrow P$, the harmonic
forms on $G/P$ go into those on $G/B_{+}$.
More precisely, 
let  $W_P \subset W$ be the subgroup of $W$ corresponding to $P$ and
let $W^P$ be the set of minimal length representatives
of the right coset space $W^P \backslash W$. For each $w \in W^P$,
Kostant constructed a harmonic form $s^{w}_{P}$, and they
have the properties stated in Kostant's Theorem 3
when $K/T$ is replaced by $K/K \cap P$ and $W$ by $W^P$.
Proposition 6.10 of \cite{ko:63} says that $\nu^* (s^{w}_{P}) =
s^w$ for $w \in W^P$. This is the third property of the Kostant
harmonic forms listed in the Introduction.
We now give a Poisson geometrical proof of this fact.

\bigskip
Let $K_P = K \cap P$. Then there is a natural identification
between $K/K_P$ and $G/P$. It is shown in \cite{lu-we:poi} that
$K_P \subset K$ is a Poisson Lie subgroup, and thus
there is an induced Poisson structure on $K/K_P$
such that the natural projection 
$\nu: K/T \rightarrow K/K_P$ is a Poisson map. This induced
Poisson structure on $K/K_P$ is also called the
Bruhat Poisson structure, because its symplectic leaves are
also exactly the Bruhat (or Schubert) cells, i.e., 
the $B_{+}$-orbits,  in $ K/K_P \cong G/P$, which are now 
indexed by elements in $W^P$. We denote this Bruhat Poisson
structure on $K/K_P$ by $\pi_{P, \infty}$.

With respect to any $K$-invariant volume form on $K/K_P$, the
modular vector field of $\pi_{P, \infty}$ is again the
infinitesimal generator of the $T$-action on $K/K_P$ in the
direction of $iH_{\rho}$ \cite{e-l-w:modular}. 
We will again use $\theta_0$
to denote this vector field. Then the same arguments we have been using
so far show that the forms $s^{w}_{P}$, for $w \in W^P$, are Poisson harmonic
with respect to the Bruhat Poisson structure $\pi_{P, \infty}$
on $K/K_P$, i.e., 
\[
(i_{\pi_{P, \infty}} d  \, - \, d i_{\pi_{P, \infty}} \, + \, i_{\theta_0})
s^{w}_{P} \, = \, 0.
\]
Because $\nu_* \pinf = \pi_{P, \infty}$ and $\nu_* \theta_0 = \theta_0$,
it follows immediately that $\nu^* s^{w}_{P}$ is also $(d, \ponf)$-harmonic.
But since $\nu^* s^{w}_{P}$ and $s^w$ have the same integrals on the 
Schubert cells in $K/T$ (Lemma 6.6 in \cite{ko:63}), and they define
the same de Rham cohomology classes in $H^{\bullet}(K/T, \C)$, they
must be equal.

\section{The $S^1$-equivariant cohomology of $K/T$}
\label{sec_equi}

In this section, we give a Poisson geometrical proof of the 
fourth property of the Kostant harmonic forms listed in the
Introduction. Namely, we show how the ring structure
on the de Rham cohomology of $K/T$ can be described using
the Kostant harmonic forms. 
We do this by using equivariant cohomology.

\bigskip
Consider the $S^1$-action on $K/T$ defined by the modular vector 
field $\theta_0$. We will use $H_{S^1}(K/T, \C)$
to denote the equivariant cohomology of $K/T$ with respect to this action.
The $S^1$-equivariant cohomology of a one point space
is the polynomial ring $\C [u]$. Thus, $H_{S^1}(K/T, \C)$
is a $\C[u]$-algebra. Since the CW-complex structure
on $K/T$ defined by the Bruhat Decomposition is $S^1$-equivariant,
we know \cite{ar:equi} that $\hus$ is  in fact a free 
$\C[u]$-module with rank equal to the
dimension of the ordinary de Rham cohomology of $K/T$, and 
correspondingly, there is a distinguished $\C[u]$-basis $\{ \sigma^{(w)}:
w \in W \}$, called the {\bf Schubert basis}, of $\hus$ that is 
characterized by the following properties: 
\begin{enumerate}
\item
$\deg \sigma^{(w)} = 2 l(w)$;
\item
$\int_{X_w} \sigma^{(w_1)} = \delta_{w, w_1}$ for $w, w_1 \in W$. Here,
$X_w$ is the closure of $\Sigma_w$ in $K/T$.
\item
Under the evaluation at $0$:
\begin{equation}
\label{eq_ev}
e_0: \, \, \hus \, \lrw \, \C \otimes_{\C[u]} H_{S^1}(K/T, \C) \, 
\cong \,  H^{\bullet}(K/T, \C)
\end{equation}
we have $\sigma^{(w)} \mapsto \sigma^w$, where $\{\sigma^w: w \in W\}$
is the Schubert basis for $H^{\bullet}(K/T, \C)$ that is dual 
to the basis of $H_{\bullet}(K/T, \C)$ defined by the closures
of the Schubert cells.
\end{enumerate}
 
\begin{rem}
\label{rem_leq}
{\em In fact, this basis has the further property that
$i_{w}^{*} \sigma^{(w_1)} = 0 \in H_{S^1}(X_w, \C)$ unless
$w_1 \leq w$, where $i_w$ is the inclusion $X_w \hookrightarrow
K/T$, and $H_{S^1}(X_w, \C)$ is the $S^1$-equivariant cohomology of
$X_w$. The partial order $\leq$ is the Bruhat order on $W$, i.e.,
$w_1 \leq w$ if and only if $\Sigma_{w_1} \subset X_w$.
}
\end{rem}

\bigskip
The $\C[u]$-algebra structure on $\hus$ is determined
by its structure constants in the Schubert basis: For
$w_1, w_2 \in W$, write
\begin{equation}
\label{eq_awww}
\cda \cdb \, = \, \sum_{w \in W} c^{w_1, w_2}_{w} \cd
\end{equation}
with $c^{w_1, w_2}_{w} \in \C[u]$. Then, since $\deg \cd 
= 2 l(w)$, we know that $c^{w_1, w_2}_{w} $ is homogeneous 
and has degree $l(w_1) + l(w_2) - l(w)$ in $u$. Moreover, since 
\[
c^{w_1, w_2}_{w} \, = \, \int_{X_w} i_{w}^{*}
\cda i_{w}^{*} \cdb,
\]
it follows from Remark \ref{rem_leq} that
\[
c^{w_1, w_2}_{w} \, = \, 0 \hspace{.3in} {\rm unless} \,
w_1 \leq w \,\,  {\rm and} \, \, w_2 \leq w.
\]
Since the map $e_0$ in (\ref{eq_ev}) is a surjective
homomorphism between algebras over $\C$, the values 
of the $c^{w_1, w_2}_{w}$'s at $u = 0$ are the structure 
constants for the algebra structure on the ordinary
de Rham cohomology $H^{\bullet}(K/T, \C)$ in the
Schubert basis $\{ \sigma^w: w \in W \}$:
\[
\sigma^{w_1} \sigma^{w_2} \, = \, \sum_{w \in W}
c^{w_1, w_2}_{w} (0) \sigma^w
\]
for $w_1, w_2 \in W$.
Thus we need to determine the polynomials $c^{w_1, w_2}_{w}$. 

\bigskip
Introduce the $\C[u]$-module
$\Hom_{\C[u]} (\hus, \, \C[u])$.
It has the obvious $\C[u]$-basis $\{\sigma_{(w)}: \, w \in W \}$
that is dual to the Schubert basis $\{ \cd: \, w \in W \}$.
On the other hand, for each $w \in W$, by regarding $w$ as 
an element in $K/T$ (via any representative of $w$ in $K$),
we have the $S^1$-equivariant map
\[
\{{\rm pt}\} \lrw K/T: \, \, {\rm pt} \Map w
\]
which induces a $\C[u]$-algebra homomorphism
\[
\psi_w: \, \, \hus \lrw \C[u],
\]
or $\psi_w \in \hls$. Write
\[
\psi_w \, = \, \sum_{w_1 \in W} d_{w_1, w} \sigma_{(w_1)}
\]
with $d_{w_1, w} \in \C[u]$. Then, by definition,
\begin{equation}
\label{eq_du}
d_{w_1, w} \, = \, \psi_w (\cda),
\end{equation}
so $d_{w_1, w}$ has degree $l(w_1)$ in $u$.
It again follows from Remark \ref{rem_leq} 
that $d_{w_1, w} = 0$ unless $w_1 \leq w$.

\begin{exam}
\label{exam_simple}
{\em
If $\alpha$ is a simple root and if
$r_{\alpha} \in W$ is the simple reflection defined by $\alpha$,
then \cite{ar:equi}
\begin{equation}
\label{eq_exam-simple}
\psi_{r_{\alpha}} \, = \, \sigma_{(1)} \, + \, 
 u\ll\alpha, \rho \gg
\sigma_{(r_{\alpha})}
\end{equation}
so
\[
d_{1, r_{\alpha}} \, =  \,1, \hspace{.3in} d_{r_{\alpha}, r_{\alpha}}
\, = \, u \ll \a, \, \rho \gg.
\]
}
\end{exam}

Notice that $\{\psi_w: w  \in W \}$ is not a $\C[u]$-basis for
$\Hom_{\C[u]} (\hus, \, \C[u])$, as is seen from (\ref{eq_exam-simple}).
Nevertheless, introduce the matrix
\begin{equation}
\label{eq_D}
D \, = \, (d_{w_1, w})_{w_1, w \in W}.
\end{equation}
Then, the localization theorem of equivariant
cohomology \cite{a-b:equi} implies that the matrix $D$ is invertible
and that the entries of the matrix $D^{-1}$ are in $\C(u)$, the algebra
of Laurent polynomials in $u$.

\bigskip
We now show how the structure
constants $c^{w_1, w_2}_{w}$ can be determined by
the matrix $D$.

\bigskip
\begin{prop}
For each $w_1 \in W$, let 
$D_{w_1}$ be the diagonal matrix with $\{d_{w_1, w}\}_{w \in W}$
as the diagonal. Let $C_{w_1}$ be the matrix whose $(w_2, w)$-entry
is $c^{w_1, w_2}_{w}$. Then
\begin{equation}
\label{eq_C-D}
C_{w_1} \, = \, D \cdot D_{w_1} \cdot D^{-1}
\end{equation}
for any $w_1 \in W$.
\end{prop}

\noindent
{\bf Proof}. Apply $\psi_w$ to both sides of
(\ref{eq_awww}). We get
\beqa
{\rm l.h.s.} & = & \psi_w( \cda \cdb) \, = \, 
 \psi_w (\cda) \psi_w (\cdb) \, = \, d_{w_1, w} d_{w_2, w}\\
{\rm r.h.s.} & = & \sum_{v \in W} d_{v, w} c^{w_1, w_2}_{v}.
\eeqa
Thus $d_{w_1, w} d_{w_2, w} = \sum_{v \in W} d_{v, w} c^{w_1, w_2}_{v}$
for all $w, w_1$ and $w_2 \in W$, which is expressed in
a more compact way in
(\ref{eq_C-D}).
\qed

\bigskip
Therefore, to determine the structure constants $\{c^{w_1, w_2}_{w}\}$, 
it is enough to determine the matrix $D$. 

\bigskip
\begin{rem}
\label{rem_D-k-k}
{\em
If we consider the $T$-equivariant cohomology $H_T(K/T, \C)$,
a $D$-matrix, whose entries are polynomials on the Cartan 
subalgebra $\fh$, can be similarly defined. This is the $D$-matrix introduced
in Section 4.21 of \cite{k-k:integral} by Kostant and Kumar. 
We refer  to it as the full $D$-matrix. 
The matrix $D$ we have here is the restriction to $uH_{\rho}$ of the
full $D$-matrix.
}
\end{rem}

\bigskip
To see how the Kostant harmonic forms $s^w, w \in W$, can be
used to calculate the matrix $D$, we recall that
these forms are also Poisson harmonic (Corollary \ref{cor_main-2})
with respect to the Bruhat Poisson structure $\pinf$. 
Thus, we can apply Theorem \ref{thm_equi} to them to get

\begin{thm}
\label{thm_main-3}
For $w \in W$, set
\begin{equation}
\label{eq_swu}
s^w(u) = i_{\exp_{\wedge} (-u \pi_{\infty})} s^{w} = s^w - u
i_{\pi_{\infty}} s^w
+ {\frac{u^2}{2!}} i_{\pi_{\infty} \wedge \pi_{\infty}} s^w + \cdots.
\end{equation}
Then

1) Each $s^w(u)$ is $S^1$-equivariantly closed;

2) The set
\begin{equation}
\label{eq_the-set}
\{{[s^w(u)] \over \lambda_w}: ~ w \in W \},
\end{equation}
 where $\lambda_w$ is the number given in (\ref{eq_lambda-w}),
 is the Schubert basis for $\hus$.
\end{thm}

\noindent
{\bf Proof.} 1) follows from Theorem \ref{thm_equi}. Clearly, 
$\deg s^w(u) = 2 l(w)$ because $\deg u = 2$. Since 
$s^w(0) = s^w$, it follows from 
Kostant's Theorem 3 in Section \ref{sec_kos-thm} that
the set in (\ref{eq_the-set}) goes to the Schubert basis
of $H^{\bullet}(K/T, \C)$ under the evaluation map $e_0$. 
We now need to show
that $\int_{X_w} s^{w_1}(u) = \delta_{w, w_1}$. By replacing
$X_w$ with a $T$-equivariant resolution of singularities $Z\to X_w$
which is an isomorphism over $\Sigma_w$, we can show
$\int_{X_w} s^{w_1}(u) = \int_{\Sigma_w} s^{w_1}(u).$
Thus, we only need to 
show  
$\int_{\Sigma_w} s^{w_1}(u) = \delta_{w, w_1}$.  This is clearly
true when $l(w_1) < l(w)$. When $l(w_1) = l(w)$, it follows again from
Kostant's Theorem 3 in Section \ref{sec_kos-thm}. It remains to 
show that $\int_{\Sigma_w} s^{w_1}(u) = 0$ when $l(w_1) > l(w)$. 
In this case, the only term in $s^{w_1}(u)$ that could possibly 
contribute to the integral is the term containing 
$i_{\pi_{\infty}^{l(w_1)-l(w)}} s^{w_1}$. But the pull-back of this
term $\tau$ to $\Sigma_w$ is zero for the following reason: 
since it is a top degree form, its value
at any point in $\Sigma_w$ is 
determined by $i_{\pi_{\infty}^{l(w)}}\tau$ since $\pi_\infty$
corresponds to a symplectic form on $\Sigma_w.$ This last expression is
zero because $\pi_{\infty}^{l(w_1)} = 0$ on 
$\Sigma_w$. Thus $\int_{\Sigma_w} s^{w_1}(u) = 0$. This shows that
the set in (\ref{eq_the-set}) is the Schubert basis for
$\hus$. 
\qed

\begin{cor}
\label{cor_D-sw}
The matrix $D$ in (\ref{eq_D})
 can be calculated from the Kostant harmonic forms
$s^w, w \in W$, by
\begin{equation}
\label{eq_D-sw}
d_{w_1, w} \, = \, {1 \over \lambda_{w_1}} u^{l(w_1)} \left(s^{w_1}, \,
{\pi_{\infty}^{l(w_1)} \over l(w_1)!} \right)(w) \, = \, 
{1 \over \lambda_{w_1}} u^{l(w_1)} (s^{w_1}, \, \, \exp_{\wedge} \pinf)(w),
\end{equation}
where $\lambda_w$ is given by (\ref{eq_lambda-w}).
\end{cor}

\noindent
{\bf Proof.} Since the $S^1$-equivariantly closed form 
${s^{w_1}(u) \over \lambda_{w_1}}$ is a
representative of the cohomology class $\cda$, the polynomial
$d_{w_1, w}$ is, by (\ref{eq_du}), the last term in the expansion
of $s^{w_1}(u)$ in (\ref{eq_swu}) evaluated at the point $w \in K/T$.
Thus we have (\ref{eq_D-sw}).
\qed

\bigskip
A description of how the entries $d_{w_1, w}$ of
the matrix $D$ can be obtained from the Kostant harmonic forms
is given in Chapter 5  (Corollary 5.6) of 
\cite{k-k:integral}. Using the formula for $\pinf (w_1)$
as given in the proof of Lemma \ref{lem_h}, it is easy
to see that our formula (\ref{eq_D-sw}) for $d_{w_1, w}$ is
the same as the one given in \cite{k-k:integral}.
We have thus given a Poisson theoretical proof of the 
fourth property of the Kostant harmonic forms listed in the
Introduction.

\bigskip
\begin{rem}
\label{rem_sara}
{\em 
S. Kumar \cite{ku:schubert} has shown that the full D-matrix can be
used to determine the singular locus of a Schubert variety  and to
determine the $T$-character of the tangent cone. We plan to
extend our work in this paper to give a Poisson construction
for the full D-matrix itself. 
This would provide a Poisson geometrical approach
to a fundamental problem in algebraic geometry. It would be interesting
to find other algebraic varieties which arise as closures of symplectic
leaves in manifolds with Poisson structures.
We remark that
S. Billey \cite{sara:d} has given an explicit formula 
for $d_{w_1, w}$ using a reduced
decomposition for $w$ as a sum over reduced 
decompositions of $w_1$ that occur as subwords of $w$.
D. Peterson \cite{dale:mit} has a much more detailed analysis of 
the structures on the $T$-equivariant cohomology $H_{T}(K/T)$.
}
\end{rem}

\bigskip
For each $w$, consider the function $F^w$ on $K/T$ given by
\begin{equation}
\label{eq_Fw}
F^w (kT) \, = \,  \left(s^{w}, \, {\pi_{\infty}^{l(w)} \over l(w)!} \right) (kT)
\, = \, (s^w, \, \, \exp_{\wedge} \pinf) (kT).
\end{equation}
Then Corollary \ref{cor_D-sw} says that 
\[
d_{w_1, w} \, = \, {1 \over \lambda_{w_1}} u^{l(w_1)} F^{w_1}(w).
\]
Thus, it is desirable to study the functions $F^w, w \in W$.
In the following, we show that they are matrix entries of
a finite dimensional representation of $K$ that is equivalent
to the Adjoint representation of $K$ on $\wedge \fg$.

\bigskip
For each $w \in W$, let $\tilde{s}^w$ be the left invariant
$2l(w)$-form on $K$ that is the pull-back of $s^w$ by the
projection $K \rightarrow K/T$.
Then $\tilde{s}^w(e)$ is a degree $2l(w)$ homogeneous element
in $\wedge \fgs$. Consider the operator ${\cal  E}$ on $\wedge \fg$
defined by
\[
{\cal E} (X) \, = \, (\exp_{\wedge} r) \wedge X \, = \, 
X + r \wedge X + {1 \over 2!} r \wedge r \wedge X + \cdots,
\]
where $r \in \fk \wedge \fk \subset \fg \wedge \fg$
is the $r$-matrix given in (\ref{eq_r}):
\[
r \, = \, {1 \over 4} \sum_{\alpha > 0} \Xa \wedge \Ya \, = \, 
{i \over 2} \sum_{\alpha > 0} \ea \wedge \eb.
\]
Then ${\cal E}$ is invertible with ${\cal E}^{-1} (X) = 
(\exp_{\wedge} (-r)) \wedge X$.
Define a representation of $K$ on 
$\wedge \fg$ by making $k \in K$ act by ${\cal E}^{-1} \circ {\rm Ad}_k
\circ {\cal E}$:
\[
k \triangleright X \, := \, \exp_{\wedge} (-r) \wedge {\rm Ad}_k
( (\exp_{\wedge} r) \wedge X) \, = \, 
\exp_{\wedge} (-r + {\rm Ad}_k r) \wedge {\rm Ad}_k X.
\]
Then we have

\bigskip
\begin{prop}
\label{prop_ad}
For any $w \in W$ and $k \in K$,
\[
F^w(kT) \, = \, (\tilde{s}^w(e), \,  \, \, \, k^{-1} \triangleright 1),
\]
where $k^{-1}\triangleright 1 =
 \exp_{\wedge} ( -r + {\rm Ad}_{k^{-1}}r).$
\end{prop}

\noindent
{\bf Proof.} Recall that $\pi = r^R - r^L$ is the Poisson structure on $K$
(see (\ref{eq_pi-on-K})) that makes $(K, \pi)$ into a Poisson Lie group,
where $r^R$ and $r^L$ are respectively the right and left invariant 
bi-vector fields on $K$ with values $r$ at the identity element $e$.
Since, by the definition of $\pinf$, the natural projection
$(K, \pi) \rightarrow (K/T, \pinf)$ is a Poisson map, we have
\beqa
F^w(kT) & = & (s^w, \, \, \exp_{\wedge} \pinf) (kT) \\
& = & (\tilde{s}^w, \, \, \exp_{\wedge} \pi)(k) \\
& = & (\tilde{s}^w, \, \, (\exp_{\wedge} r)^R (\exp_{\wedge} (-r))^L) \\
& = & (\tilde{s}^w(e), \, \, (\exp_{\wedge}(-r)) \wedge {\rm Ad}_{k^{-1}}
(\exp_{\wedge} r) )\\
& = & (\tilde{s}^w(e), \, \, k^{-1} \triangleright 1).
\eeqa
\qed

\bigskip
In a separate paper, we plan to study further properties of these
functions $F^w$ for $w \in W$, as well as similar functions 
for the full $D$-matrix.

\section{Kostant's harmonic forms as limits of Hodge harmonic forms}
\label{sec_family}

In this section, we first introduce a family 
of Poisson (in fact symplectic) structures $\pl$ on $K/T$, where
$\lambda \in \fa$ is regular. We study Poisson harmonic forms
for each $\pl$ and show that they are the same as Hodge harmonic forms
for certain Hermitian metrics $h_{\lambda}$ on $K/T$. 
We show that the Kostant operator $\partial$ is
the limit as $\lambda \rightarrow \infty$ of the 
adjoint operators $d_{\ast, \lambda}$ of $d$ with
respect to the Hermitian metrics $h_{\lambda}$.
As a consequence, the
Laplacian $S = d \partial + \partial d$ defined by Kostant is the limit of the 
Hodge Laplacian of $d$ with respect to the Hermitian metrics
$h_{\lambda}$.
This will enable us to give a new proof of Kostant's
theorem that $d$ and $\partial$ are disjoint.
Finally, we define Hodge harmonic forms $s^{w}_{\lambda}$
for $w \in W$  and 
show that they tend to the Kostant harmonic forms $s^w$ as $\lambda 
\rightarrow \infty$.

\subsection{The family of symplectic structures $\pl$}
\label{sec_pi-lambda}
We first recall some facts about the Poisson Lie group
$(K, \pi)$. Details can be found in \cite{lu-we:poi}.

\bigskip
Let $G = K_{\Bbb C}$ be the complexification of $K$, but here regarded
as a real Lie group. Let $A$ and $N$ be the connected Lie subgroups
of $G$ with Lie algebras $\fa$ and $\fn$ respectively. Then $G = KAN$ 
is an Iwasawa decomposition of $G$ as a real semi-simple Lie group.
The group $AN$ has a unique Poisson structure $\pi_{AN}$ 
making $(AN, \pi_{AN})$ into the 
dual Poisson Lie group of $(K, \pi)$. The group $G$ is then
the double group for $(K, \pi)$. The decomposition $G = KAN$
gives rise to the left action of $K$ on $AN$:
\[
K \times AN \lrw AN: \, \, (k, b) \Map k \cdot b :=b_1, \hspace{.3in}
{\rm if}\,  \, \, b k^{-1} = k_1 b_1 \,\, {\rm for}\, \, k_1 \in K
\,\, 
{\rm and}\, \, b_1 \in AN.
\]
It is called the (left) dressing action of $K$ on $AN$, and its
orbits are exactly
the symplectic leaves of $\pi_{AN}$ in $AN$. 
Thus each dressing orbit inherits a symplectic, and thus Poisson,
structure as a symplectic leaf. Since the dressing action
is Poisson \cite{sts:dressing} \cite{lu-we:poi}, these
dressing orbits are examples of
$(K, \pi)$-homogeneous Poisson spaces.
Let $\lambda \in \fa$ be regular and consider the element
$\el \in A$.
The stabilizer subgroup of $K$ in $AN$
at $\el$ is $T$. Thus, by identifying $K/T$ with the dressing orbit
through $\el$, we get a Poisson structure on $K/T$ which is in fact
symplectic.

\begin{nota}
\label{nota_symplectic}
{\em
For $\lambda \in \fa$ regular, we will use
$\pi_{\lambda}$ to denote
the Poisson structure on $K/T$ obtained by identifying
$K/T$ with the symplectic leaf in $AN$ through the point $\el$.
We call it the
{\it dressing orbit Poisson structure corresponding to $\el$}.
}
\end{nota}

The following proposition is proved in \cite{lu:cdyb}.

\begin{prop}
\label{prop_pl-explicit}
The Poisson structure $\pl$ on $K/T$ is explicitly given by
\[
\pl \, = \, \left( \sum_{\alpha > 0} {1 \over 1-e^{2\alpha (\lambda)}}
{\Xa \wedge \Ya \over 2} \right)^{0} \, + \, \pinf,
\]
where the first term is the $K$-invariant bi-vector field on $K/T$
whose value at $\eo = eT$ is the expression given in the parentheses.
The modular vector field of $\pl$ with respect to
a $K$-invariant volume form on $K/T$ is again the vector
field $\theta_0$ given in (\ref{eq_theta-0}).
\end{prop}

\begin{nota}
\label{nota_pl}
{\em
We use $\pal$ to denote the Koszul-Brylinski
operator defined by the Poisson structure $\pl$ and $\mu_0$:
\[
\pal \, = \, \partial_{\pi_{\lambda}, \mu_0} \, = \,
i_{\pi_{\lambda}} d \, - \, d i_{\pi_{\lambda}} \, + \,
i_{\theta_0}:\, \,
\Omega^q(K/T) \lrw \Omega^{q-1}(K/T).
\]
}
\end{nota}

\bigskip
The following proposition explains the notation $\pinf$ and
$\ponf$ we have 
given to the Bruhat Poisson structure and its
associated Koszul-Brylinski operator:

\begin{prop}
\label{prop_limit-1}
For any $\lambda \in \fa$,
\[
\lim_{t \rightarrow + \infty} \pi_{\lambda + tH_{\rho}} \, = \, \pinf
\hspace{.3in} {\rm and} \hspace{.3in}
\lim_{t \rightarrow + \infty} \partial_{\lambda + tH_{\rho}} \, = \, 
\ponf
\]
in the $C^{\infty}$-topology on the space of tensors on $K/T$.
\end{prop}

\noindent
{\bf Proof.} This is because $\rho(\alpha) > 0$ for each $
\alpha > 0$, and thus 
\[
\lim_{t \rightarrow + \infty} {1 \over 1 \, - \, 
e^{2 \alpha (\lambda + t H_{\rho})}} \, = \, 0.
\]
\qed

\bigskip
As in the case of the Bruhat Poisson structure $\pinf$, 
since $\pl$ is $(K, \pi)$-homogeneous 
and since the volume form $\mu_0$ is $K$-invariant,
the operator $\pal$ leaves invariant the space 
$\Omega(K/T)^K$ of $K$-invariant (real) differential forms
on $K/T$. Therefore, we have the mixed complex
$(\Omega(K/T)^K, \, d, \, \pal)$ for each regular 
element $\lambda \in \fa$.

\bigskip
\begin{rem}
\label{rem_simple-2}
{\em
The restriction of $\pal$ to the space
$\Omega(K/T)^K$ of $K$-invariant forms on $K/T$ has
a simpler expression, just like
the case with $\ponf$ (see Remark \ref{rem_simple-1}).
Namely, consider the element
\[
r_{\lambda} \, = \,{1 \over 4}  \sum_{\alpha > 0}
{e^{\alpha(\lambda)} + e^{-\alpha(\lambda)} \over
e^{\alpha(\lambda)} - e^{-\alpha(\lambda)}} \Xa \wedge \Ya
\]
as an element in $\wedge^2 (\fk / \ft)$ and extend
it $K$-invariantly to a bi-vector field on $K/T$, which will
be denoted by the same letter.
Then 
\[
\pal \, = \, - i_{r_{\lambda}} d \, + \, 
d i_{r_{\lambda}}
\]
on $\Omega(K/T)^K$. Notice that this
is not true on the space $\Omega(K/T)$
of all differential forms on $K/T$.
The element $r_{\lambda}$ is the skew-symmetric part of 
a {\it classical dynamical $r$-matrix} studied by
Etingof and Varchenko in \cite{e-v:cdyb}. Connections
between such $r$-matrices and homogeneous Poisson structures
are given in \cite{lu:cdyb}.
}
\end{rem}

\bigskip
As in the case for the Bruhat Poisson structure $\pinf$,
we will first identify the operator 
$\pal$ on $\Omega(K/T)^K$ with the Chevalley-Eilenberg
boundary operator for some Lie algebra. For this 
purpose, we introduce the following real Lie subalgebra
of $\fg$:
\[
\fll \, = \, {\rm Ad}_{e^{\lambda}} \fk \, \in \, \fg.
\]
It is isotropic with respect to the 
bi-linear form $2{\rm Im} \ll \, , \, \gg$ and is maximal with
this property since it  has half of the 
real dimension of $\fg$.
Moreover, $\fll \cap \fk = \ft$. Under the correspondence
established by Drinfeld \cite{dr:homog}
 between maximal isotropic Lie algebras
$\fl$ of $\fg$ satisfying $\fl \cap \fk = \ft$ and $(K, \pi)$-homogeneous
Poisson structures on $K/T$, the Lie algebra $\fll$ corresponds
to the Poisson structure $\pl$. This fact is proved in
\cite{lu:homog} and \cite{lu:cdyb}. 

\begin{rem}
\label{rem_fl-limit}
{\em
In the case of $\pinf$, the corresponding 
maximal isotropic Lie subalgebra of $\fg$ is $\ft + \fn$.
Notice that
\beqa
{\rm Ad}_{e^{\lambda}} \Xa &  = &  
e^{\alpha (\lambda)} (\ea \, - \, e^{-2 \alpha (\lambda)} \eb) \\
{\rm Ad}_{e^{\lambda}} \Ya &  = &   
e^{\alpha (\lambda)}  ( i\ea \, + \, e^{-2 \alpha (\lambda)} i\eb).
\eeqa
Thus
\[
\lim_{t \rightarrow + \infty} \fl_{\lambda+ tH_{\rho}} \, = \, \ft + \fn
\]
in the Grassmannian of maximal isotropic 
Lie subalgebras of $\fg$.
This fact corresponds to the statement in Proposition 
\ref{prop_limit-1}}.
\end{rem}

\bigskip
As always, we identify
\[
\te \, \cong \, \fk / \ft \hspace{.3in}
{\rm and} \hspace{.3in} 
\tse \, \cong \, (\fk / \ft)^{*}.
\]
The bi-linear form $2{\rm Im} \ll \, , \, \gg$ now gives an identification
\begin{equation}
\label{eq_Il}
\il: \, \, \fll / \ft \, \stackrel{\sim}{\lrw} \, \tse: \, \,
(\il (x + \ft), \, y + \ft) \, = \, 
2{\rm Im} \ll x, \, y \, \gg,
\end{equation}
where $x + \ft \in \fll / \ft$ and $y + \ft \in \fk / \ft$.
Its dual map
\begin{equation}
\label{eq_ils}
\ils: \, \, \te \, \stackrel{\sim}{\lrw} \, 
(\fll / \ft)^*
\end{equation}
is defined by the same bi-linear form. They give rise to identifications,
still denoted by the same letters,
\beqa
& & \il: \, \, (\wedge^q (\fll / \ft))^T \, \stackrel{\sim}{\lrw} \,
\Omega^q(K/T)^K \\
& & \ils: \, \, \chi^q(K/T)^K \, \stackrel{\sim}{\lrw} \, (
\wedge^q (\fk / \ft)^*)^T.
\eeqa
Denote by
\[
d_{{\frak l}_{\lambda}}: \, \, \wedge^q \fl_{\lambda}^{*} \lrw 
\wedge^{q+1} \fl_{\lambda}^{*}
\]
the Chevalley-Eilenberg coboundary operator for the Lie algebra
$\fll$, and by the same letter its restriction
to the subspace $(\wedge (\fk / \ft)^*)^T \subset \wedge \fl_{\lambda}^{*}$.
Denote by
\[
\bfl: \, \, \wedge^q \fll \lrw \wedge^{q-1} \fll
\]
the Chevalley-Eilenberg boundary operator for the Lie algebra
$\fll$, and by the same letter the induced operator
on $(\wedge (\fll / \ft))^T$ that calculates the relative Lie
algebra homology of $\fll$ relative to $T$ \cite{kn:cohom}. 
The following statement  is parallel to that in
Theorem \ref{thm_ponf-n}. 

\begin{prop}
\label{prop_pl-fll}
We have
\begin{equation}
\label{eq_lambda-d}
(\ils)^{-1} \dfl \ils \, = \, \delta_{\pi_{\lambda}} \, = \, 
[\pl, \, \bullet]
\end{equation}
as degree $1$ operators on $\chi(K/T)^K$, and
\begin{equation}
\label{eq_lambda-partial}
\il \bfl I_{\lambda}^{-1} \, = \, \pal
\end{equation}
as degree $-1$ operators on $\Omega(K/T)^K$.
\end{prop}

\noindent
{\bf Proof.} As with Theorem \ref{thm_ponf-n},
this statement is a special case of a general fact about
Poisson homogeneous spaces \cite{lu:homog} \cite{lu:cdyb}. 
A direct proof similar to that for
Theorem \ref{thm_ponf-n} can be given. We omit it here.
\qed

\bigskip
Unlike $\ft + \fn$, the Lie algebra $\fll$
is isomorphic to the Lie algebra $\fk$. As a result, 
the operator $\pal$ is conjugate to the operator $d$ on
$\Omega(K/T)^K$ as is shown now. Introduce
\begin{equation}
\label{eq_m-lambda}
m_{\lambda} \, = \, \sum_{\alpha > 0}
{e^{\alpha (\lambda)}  \over 2(1 -  e^{2\alpha(\lambda)})}
 \, \Xa \wedge \Ya,
\end{equation}
considered as an element in $\wedge^2 (\fk /\ft) \cong
\wedge^2(T_{e}(K/T))$. See Remark \ref{rem_m-construction}
for the construction of $m_{\lambda}$ from $\pl$.
 Since $m_{\lambda}$ is non-degenerate,
the map
\[
\tilde{m}_{\lambda}: \, \, T^{*}_{e}(K/T) \lrw 
T_{e}(K/T): \, \,
\xi \Map \xi \backl \, m_{\lambda}
\]
is a vector space isomorphism. Extend it to
\[
\tilde{m}_{\lambda}: \, \, \wedge^q (T^{*}_{e}(K/T))
\lrw \wedge^q \te: \, \,
\tml(\xi_1 \wedge \cdots \wedge \xi_q) \, = \,
\tml(\xi_1) \wedge \cdots \wedge \tml(\xi_q).
\]
Recall that $\mu_0$ is a fixed $K$-invariant volume
form on $K/T$. Regard $\mu_0$ as in $\wedge^n(\tse)$,
where $n = \dim_{\Bbb R} K/T$.
Introduce the following analog of the Hodge $\ast$-operator
(see \cite{by:poi}):
\begin{equation}
\label{eq_star-M}
\star_{\lambda}: \, \, \wedge^q(\tse) \lrw \wedge^{n-q}(\tse): \, \,
\star_{\lambda} (\xi) \, = \, \tml(\xi) \backl \, \mu_0.
\end{equation}
It is easy to show that
\[
\star^{2}_{\lambda} \, = \, \left({m_{\lambda}^{{n \over 2}} 
\over ({n \over 2}!)}, 
\,\, \mu_0 \right)^2
{\rm id}.
\]
Since $\star_{\lambda}$ is $T$-invariant, it can be regarded as
an operator on $\Omega(K/T)^K \cong \wedge (\tse)^T$.
 
\bigskip
\begin{prop}
\label{prop_partial-d-star}
We have
\[
\partial_{\lambda} \, = \, (-1)^q \star_{\lambda}^{-1} d \, 
\star_{\lambda}:\, \,
\Omega^q(K/T)^K \lrw \Omega^{q-1}(K/T)^K.
\]
\end{prop}

\bigskip
\noindent
{\bf Proof.}
As always, we identify $\te$ with $\fk / \ft$, so we can regard
both $d$ and $\partial_{\lambda}$ as operators on
$(\wedge (\fk / \ft)^*)^T$. As such, the operator $d$ is nothing but the
restriction to $(\wedge (\fk / \ft)^*)^T \subset \wedge \fk^*$
of the Chevalley-Eilenberg coboundary operator for the
Lie algebra $\fk$.
 
Denote by
\[
b_{\frak k}: \, \, \wedge^q(\fk / \ft)^T \lrw \wedge^{q-1}
(\fk / \ft)^T
\]
the Chevalley-Eilenberg boundary operator that calculates the relative
Lie algebra homology of $\fk$ relative to $\ft$ (see \cite{kn:cohom}).
Then, since ${\rm Ad}_{e^{\lambda}}: \fk \rightarrow \fl_{\lambda}$
is a Lie algebra isomorphism, we know from Proposition
\ref{prop_pl-fll} that
\[
\partial_{\lambda} \, = \, (I_{\lambda} \circ {\rm Ad}_{e^{\lambda}})
\circ b_{\frak k} \circ (I_{\lambda} \circ {\rm Ad}_{e^{\lambda}})^{-1}.
\]
Now it is easy to check that
\[
(I_{\lambda} \circ {\rm Ad}_{e^{\lambda}})^{-1}: \, \,
(\fk / \ft)^* \lrw \fk / \ft
\]
is given by
\beqa
X^{\alpha} & \Map &{ e^{\alpha (\lambda)} \over  
2(1 - e^{2\alpha (\lambda)})}
 \,\Ya \\
Y^{\alpha} & \Map &-{ e^{\alpha (\lambda)} \over  
2(1 - e^{2\alpha (\lambda)})}
 \, \Xa,
\eeqa
where $\{X^{\alpha}, \, Y^{\alpha}: \, \alpha > 0 \}$
is the basis of $(\fk / \ft)^*$ dual to the basis
$\{\Xa, \, \Ya: \, \a  > 0 \}$ of $\fk / \ft$.
Thus,
\[
(I_{\lambda} \circ {\rm Ad}_{e^{\lambda}})^{-1} \, = \, \tml: \, \,
(\fk / \ft)^* \lrw \fk / \ft.
\]
Hence
\[
\partial_{\lambda} \, = \, \tml^{-1} \circ b_{\frak k} \circ \tml: \, \,
(\wedge^q(\fk / \ft)^*)^T \lrw (\wedge^{q-1}(\fk / \ft)^*)^T,
\]
or
\[
\tml (\partial_{\lambda} (\xi)) \, = \, b_{\frak k}(\tml(\xi)) \, \in
\, (\wedge^{q-1}(\fk / \ft))^T
\]
for $\xi \in (\wedge^q(\fk / \ft)^*)^T$.
But
\beqa
\tml (\partial_{\lambda} (\xi)) \backl \, \mu_0 & = &
\star_{\lambda} \partial_{\lambda} (\xi), \\
b_{\frak k}(\tml(\xi)) \backl \, \mu_0 & = & (-1)^q d(\tml(\xi) \backl \,
\mu_0) \, = \, (-1)^q d \star_{\lambda} (\xi).
\eeqa
Therefore,
\[
\star_{\lambda} \partial_{\lambda} (\xi) \, = \,
(-1)^q d \star_{\lambda} (\xi).
\]
\qed
 
\begin{rem}
\label{rem_with-B}
{\em
Our result in Proposition \ref{prop_partial-d-star} should be 
compared with Theorem 2.2.1 in \cite{by:poi}, where the operator
$i_{\pi_{\lambda}} d - d i_{\pi_{\lambda}} = \ponf - i_{\theta_0}$
is related to $d$ by a $\ast$-operator defined by the symplectic structure
$\pl$.
}
\end{rem}

\begin{rem}
\label{rem_m-construction}
{\em
At this point, we want to explain how $m_{\lambda}$
is related to the Poisson structure $\pl$ or
more precisely $\pl(e) \in
\wedge^2 \te$. The element $e^{\lambda}\in AN$
acts on $\fn = \fn_{+}$ by the Adjoint action. By identifying
$\tse$ with $\fn$ using the map $I$ in (\ref{eq_I}), we get
an action of $e^{\lambda}$ on $\tse$. Then for $\xi, \eta \in 
\tse$, we have
\[
m_{\lambda} (\xi, \eta) \, = \, 
\pl(e) ({\rm Ad}_{e^{\lambda}} \xi, \, \eta)
\, = \, \pl(e) (\xi, \, {\rm Ad}_{e^{\lambda}} 
\eta).
\]
}
\end{rem}

\bigskip
Recall now that $J$ is the complex structure
on $\te \cong \fk / \ft$ coming from the identification
$\fk / \ft \cong \fg / \fb_{-}$. The same letter also denotes
the dual complex structure on $\tse \cong (\fk / \ft)^*$.
Introduce the following bi-linear form on $\tse$:
\[
g_{\lambda}(\xi, \, \eta) \, = \,
m_{\lambda} ( \xi, \, J\eta)  \, = \, - m_{\lambda}(J \xi, \,  \eta).
\]
It is symmetric, and
\[
g_{\lambda}(X^{\alpha}, \, X^{\alpha}) \, = \,
g_{\lambda}(Y^{\alpha}, \, Y^{\alpha}) \, = \, 
{e^{\alpha (\lambda)} \over 2(e^{2 \alpha (\lambda)} - 1)}
\]
and $0$ otherwise. Thus, when $\lambda \in \fa$ is dominant,
the bilinear form $g_{\lambda}$ on $\tse$ is symmetric and positive
definite.
 
\bigskip
Without loss of generality, we assume
that $\lambda$ is dominant from now on. Extend $g_{\lambda}$
$K$-invariantly to all of $T^*(K/T)$. Its inverse
on $T(K/T)$ is then a Riemannian metric on $K/T$. We use the
same letter $g_{\lambda}$ to denote the 
multi-linear extension of $g_{\lambda}$
to $\wedge T^* (K/T)$. 

\bigskip
Consider now the Hodge $\ast$-operator
associated to $g_{\lambda}$. It is defined by
\[
\xi \wedge \ast_{\lambda} \eta \, = \, \eta \wedge \ast_{\lambda} 
\xi \, = \, 
g_{\lambda}(\xi, \eta) \mu_0,
\hspace{.2in} {\rm for} \hspace{.1in} 
\xi, \eta \in \Omega^q(K/T).
\]
It is also given by
\[
\ast_{\lambda} \xi \, = \, \tilde{g}_{\lambda}(\xi) \backl \mu_0
\]
where $\tilde{g}_{\lambda}(\xi)$ is defined by
$\tilde{g}_{\lambda}(\xi) (\eta) = g_{\lambda} (\xi, \eta)$
for $\xi, \eta \in \Omega^q (K/T)$. We have
\[
\ast_{\lambda}^{2} \xi \, = \, (-1)^q g_{\lambda}(\mu_0, \,
\mu_0) \xi, \hspace{.2in} {\rm for} \hspace{.1in}  
\xi  \in \Omega^q(K/T).
\]
Since $g_{\lambda}$ is $T$-invariant, 
we can regard $\ast_{\lambda}$ as an operator on $\Omega(K/T)^K$:
\[
\ast_{\lambda}: \, \, \Omega^q(K/T)^K \lrw
\Omega^{n-q}(K/T)^K.
\]

Consider now the symmetric bilinear form $( \, \, \, )_{\lambda}$
on $\Omega(K/T)$ defined by
\[
(\xi, \, \eta) \, = \, \int_{K/T} g_{\lambda} (\xi, \, \eta) \mu_0.
\]
It is positive definite. Denote by $d_{\ast,\lambda}$
 the adjoint
operator of $d$ with respect to $( \, \, \, )_{\lambda}$. It
is standard to show that 
\[
d_{\ast, \lambda} \, = \, (-1)^q \ast_{\lambda}^{-1} d \ast_{\lambda}: \, \,
\Omega^{q}(K/T) \lrw \Omega^{q-1}(K/T).
\]

\begin{prop}
\label{prop_pl-dt}
The Koszul-Brylinski operator $\pal$ and the adjoint operator 
$d_{\ast,\lambda}$ of $d$ are related via the complex structure $J$ by
\[
d_{\ast,\lambda} \, = \, J \pal J^{-1}.
\]
\end{prop}

\noindent
{\bf Proof.} It follows from the definition of $g_{\lambda}$ that
\[
\ast_{\lambda} \, = \, \star_{\lambda} J^{-1}.
\]
Thus, from Proposition \ref{prop_partial-d-star} we know
that, as operators
from $\Omega^q(K/T)^K$ to $
\Omega^{q-1}(K/T)^K$,
\[
J \pal J^{-1} \, = \, 
(-1)^q J \star_{\lambda}^{-1} d \star_{\lambda} J^{-1}
\, = \, (-1)^q \ast_{\lambda}^{-1} d \ast_{\lambda} 
\, = \, d_{\ast,\lambda}.
\]
\qed

We now look at the complex picture. We complex linearly extend
both $m_{\lambda}$ and $g_{\lambda}$ to bi-linear forms on
the complexified cotangent bundle $T_{\Bbb C}^{*} (K/T)$. All the operators:
$\star_{\lambda}, \ast_{\lambda}, J, \pl, d$ and $d_{\ast, \lambda}$
are also complex linearly extended to operators on
$\Omega(K/T, \C)$ but will still be denoted
by the same letters. Let
$h_{\lambda}$ ($H$ for Hermitian)
 be the Hermitian extension of $g_{\lambda}$ to
$T_{\Bbb C}^{*} (K/T)$, i.e.,
\[
h_{\lambda} (\xi_1, \, \xi_{2}) \, = \, g_{\lambda}
(\xi_1, \, \bar{\xi}_2),
\]
where $-$ is the complex conjugation on $T_{\Bbb C}^{*} (K/T)$. 
Correspondingly, we
have the Hermitian inner product on $\Omega(K/T, \C)$:
\begin{equation}
\label{eq_hermitian}
\la\xi, \, \eta \ra_{\lambda}  \, = \, \int_{K/T} h_{\lambda}
(\xi, \, \eta) \mu_0 
\, = \,  \int_{K/T} g_{\lambda}
(\xi, \, \bar{\eta}) \mu_0.
\end{equation}
Then the operator $d_{\ast,\lambda}$ is also adjoint to $d$ with respect to
the Hermitian inner product $\la \, \, , \, \, \ra_{\lambda}$.
Let
\begin{equation}
\label{eq_S-lambda}
S_{\lambda} \, = \, d d_{\ast, \lambda} \, + \, d_{\ast, \lambda} d: 
\, \, \Omega^q(K/T, \C) \lrw \Omega^q(K/T, \C)
\end{equation}
be the (Hodge) Laplacian of $d$. It is then self-adjoint and non-negative
definite with respect to $\la \, \, , \, \, \ra_{\lambda}$. 

\bigskip
Recall from Section \ref{sec_partial} that
$\partial$ is the Kostant operator on $C = \Omega(K/T, \C)^K$
of degree $-1$ and that $S = d \partial + \partial d$
is called the ``Laplacian" of $d$ and $\partial$ in \cite{ko:63}.
Our next theorem says that 
$\partial$ is the limit of the operators $d_{\ast, \lambda}$
and that Kostant's Laplacian  $S$ is a limit of
the Hodge Laplacians $S_{\lambda}$.

\begin{thm}
\label{thm_limit-S}
With respect to the vector space topology on $C = \Omega(K/T, \C)^K$, we
have
\beqa
& & \lim_{t \rightarrow +\infty} d_{\ast,\lambda + t H_{\rho}} \, = \, 
\partial;
\\
& & \lim_{t \rightarrow +\infty} S_{\lambda + t H_{\rho}} \, = \, S.
\eeqa
\end{thm}

\noindent
{\bf Proof.} This follows directly from Propositions \ref{prop_pl-dt}
and 
\ref{prop_limit-1} 
and Theorem \ref{thm_main-1}.
\qed

It is easy to see that the Hermitian
metric $h_{\lambda}$ on $(\tse)_{\Bbb C}$
becomes the following one on $\fn_{-} \oplus \fn_{+}$
under the identification $I$:
\[
(h_{\lambda} \circ I) (\ea, \, \ea) \, = \,
(h_{\lambda} \circ I) (\eb, \, \eb) \, = \,
{e^{\alpha(\lambda)} \over e^{2\alpha (\lambda)} - 1},
\]
and $0$ otherwise. Compare $h_{\lambda}$ with the Hermitian
metric $h$ on $(\tse)_{\Bbb C}$ defined in Section 
\ref{sec_kos-thm} in (\ref{eq_h-I}). Note that
\[
\lim_{t \rightarrow +\infty} h_{\lambda + t H_{\rho}} \, = \, 0.
\]

\bigskip 
In Section \ref{sec_poi-hodge}, we will use this explicit 
formula for $h_{\lambda}$ to give another proof of Theorem 
\ref{thm_limit-S}.

\subsection{Another proof of the disjointness of $d$ and $\partial$}
\label{sec_another-proof}

Recall that Theorem 4.5 in \cite{ko:63} says that
the two operators $d$ and $\partial$ on
$C = \Omega(K/T, \C)^K$ are disjoint in the sense
that $\Im (d) \cap \Ker (\partial)
= \Im (\partial) \cap \Ker (d) = 0$. Note that the operator
$d$ is always disjoint from its adjoint operator
with respect to any Hermitian inner product. Although
$\partial$ is in general not the adjoint operator of $d$ with respect to
any Hermitian metric, as we have seen from Example \ref{exam_not-adjoint},
Theorem \ref{thm_limit-S} says that it 
is the limit of such operators 
with respect to the family of Hermitian metrics
$h_{\lambda}$. In this section, we make use of
this fact and some simple linear algebra arguments to give another
proof of the disjointness of $d$ and $\partial$.

\bigskip
Consider again the operator $S = d \partial + \partial d$ on $C$.
Recall that $S = L + E$, where
$L = \partial \delta + \delta \partial$ is the
Hodge Laplacian of $\partial$ with respect to the
Hermitian inner product $h$ on $C$ defined
in Section \ref{sec_kos-thm} (see (\ref{eq_h-I})),
and that the operator $E_I := I^{-1}EI$
on $\wedge (\fn_{-} \oplus \fn_{+})^T$ has the explicit formula
given in (\ref{eq_E}). 

\bigskip
We need two lemmas from linear algebra.

\begin{lem}
\label{lem_linear-1}
Let $C$ be a finite dimensional complex vector space. Assume
that $L$ and $E$ are two linear operators on $C$ such that
there exists a basis $\{e_1, e_2, ..., e_m \}$ of $C$ in which
$L$ is diagonal and $E$ strictly lower triangular. Let $L_1$ be the 
Green's operator of $L$, i.e., 
$L_1 (e_j) = 0$ if $L(e_j) = 0$ and $L_1(e_j) = {1 \over \lambda_j}
e_j$ if $L(e_j) = \lambda_j e_j$ and $\lambda_j \neq 0$, 
for $j = 1, ..., m$. 
Set $S = L + E$. Then, we have two
injective linear maps
\beqa
\phi: & &  \Ker (S) \lrw \Ker (L): \, \, x \Map (1 + L_1 E) x;\\
\psi: & & \Im (L) \lrw \Im (S): \, \, x \Map (1 + E L_1) x.
\eeqa
If, in addition, $\phi$ is an isomorphism, then so is $\psi$, and
$\Ker (S) \cap \Im (S) = 0$.
\end{lem}

\noindent
{\bf Proof.} For $x \in \Ker (S)$, let $y = (1 + L_1 E) x$. Using
$S(x) = 0$, we have $L(y) = L(x) + LL_1 E (x) = -E(x) + 
LL_1 E (x)  \in \Ker (L)$.
But $\Ker (L) \cap \Im (L) = 0$ since $L$ is semi-simple. Thus 
$L(y) = 0$. Now for $x \in \Im (L)$, set $z = (1 + E L_1) x$.
Then using $x = LL_1 (x)$, we get
$z =  LL_1(x) + EL_1 (x) = SL_1 (x) \in \Im (S)$.
Since $L_1 E$ and $E L_1$ are both nilpotent, the operators
$1 + L_1 E$ and $1 + EL_1$ on $C$ are both invertible. Hence
both $\phi$ and $\psi$ are injective.

If $\phi$ is an isomorphism, then, by reason of dimensions, 
$\psi$ is also an isomorphism. Suppose now $x \in \Ker (S) \cap
\Im (S)$. We want to show that $x = 0$. Since
$x \in \Ker (S)$, we have $\phi(x) = (1 + L_1 E)x \in \Ker (L)$.
Since $x \in \Im (S)$ and since $\psi$ is an isomorphism, there
exists $y \in \Im (L)$ such that $x = \psi (y) = (1 + EL_1)y$.
Thus
\begin{equation}
\label{eq_phi-x}
\phi(x) \, = \, (1 + L_1 E) (1 + EL_1)y \, = \, y + L_1 E(y) 
+ EL_1 (y) + L_1 E^2 L_1(y).
\end{equation}
Suppose that $x \neq 0$. Then $y \neq 0$.
Write $y = y_1 e_1 + y_2 e_2 + \cdots + y_m e_m$, and
let $ z_0 := y_{i_0} e_{i_0}$ be the first non-zero term in this 
expression.
Since $y \in \Im (L)$ and since $L$ is diagonal in the $\{e_j \}$ 
basis, we know that $z_0 \in \Im (L)$. But since $E$ is strictly lower
triangular in the $\{e_j \}$ basis, we know from (\ref{eq_phi-x})
that $\phi(x) = z_0 + z_1$ for some $z_1 \in {\rm span} \{e_j: 
j > i_0\}$. Since $\phi(x) \in \Ker (L)$, we must have $z_0 \in 
\Ker (L)$. Thus $z_0 = 0$ because $\Ker(L) \cap \Im(L) = 0$.
This is a contradiction. Hence $x = 0$, and $\Ker (S) \cap \Im (S)
= 0$.
\qed
 
The proof of the following lemma is trivial.

\begin{lem}
\label{lem_rank}
Suppose that $V_1$ and $V_2$ are two  finite dimensional vector spaces
and that $A_t \in \Hom(V_1, V_2)$, for $t > 0$, is a smooth family of
linear operators from $V_1$ to $V_2$ such that
\[
\lim_{t \rightarrow +\infty} A_t \, = \, A \, \in \, \Hom(V_1, V_2).
\]
If $\dim (\Ker (A_t)) = n_0$ for some integer $n_0$ for all $t > 0$,
then $\dim (\Ker (A)) \geq n_0$.
\end{lem}

\bigskip
Lemma \ref{lem_linear-1} can be applied to 
the operators $S, L$ and $E$ on $C =\Omega(K/T, \C)^K$ since
all the assumptions are satisfied. Thus
we have the linear maps $\phi: \Ker (S)
\rightarrow \Ker (L)$ and $\psi: \Im (L) \rightarrow \Im(S)$
which are both injective.
The Green's operator $L_1$ for $L$, when considered as
an operator on $\wedge (\fn_{-} \oplus \fn_{+})^T$ via its
identification with $C$ by $I$, is the operator
$L_0$ given in Section \ref{sec_kos-thm}.
We use $m_0$ to denote the common number
\[
m_0 \, = \, \dim (\Ker (L)) \, = \, \dim H(C, \partial) \, = \, 
|W| \, = \, \dim (H^{\bullet}(K/T, \C)),
\]
where $|W|$ is the number of elements in the Weyl group $W$.

\bigskip
\noindent
{\bf Claim 1.} We have $\dim (\Ker (S)) = m_0$. Consequently,
the map  $\phi: \Ker (S) \rightarrow \Ker (L)$ is
an isomorphism, and  $ \Ker (S) \cap \Im (S) = 0$.

\bigskip
\noindent
{\bf Proof.}
Consider the family of operators $S_t := S_{tH_{\rho}} 
\in \End(C)$. By Theorem \ref{thm_limit-S}, we have
\[
\lim_{t \rightarrow +\infty} S_t \, = \, S.
\]
But we know from the usual Hodge theory that $\dim (\Ker (S_t)) =
m_0$ for all $t > 0$. Thus by Lemma \ref{lem_rank}, we know
that $\dim (\Ker (S)) \geq m_0$. We also know that
$\dim (\Ker (S)) \leq m_0$ because $\phi: \Ker (S) \rightarrow
\Ker (L)$ is injective. Thus $\dim (\Ker (S)) = m_0$. 
It follows from Lemma \ref{lem_linear-1} that
$\phi$ is an isomorphism and that
$\Ker (S) \cap \Im (S)  =  0$.
\qed

\bigskip
\noindent
{\bf Claim 2.} We have  $\Im (S) = \Im (d) + \Im (\partial)$.

\bigskip
\noindent
{\bf Proof.}
Clearly, $\Im (S) \subset \Im (d) + \Im (\partial)$.
Let $m = \dim (C)$. Since $\dim (\Ker (S)) = m_0$,
we have 
\[
\dim (\Im (S)) \, =  \, m - m_0  \, = \, 
m - \dim (\Ker (d)) + \dim (\Im (d))  = 2 \dim (\Im (d)).
\]
On the other hand, 
\[
\dim (\Im (d) \, + \, \Im (\partial)) \, \leq \, 
\dim (\Im (d)) \, + \,  \dim(\Im (\partial)).
\]
Consider now the family of operators $d_{\ast, tH_{\rho}}$ on $C$ 
for $t > 0$. By 
Theorem \ref{thm_limit-S}, 
\[
\lim_{t \rightarrow +\infty} d_{\ast, tH_{\rho}} \, = \, 
\partial.
\]
Thus we know from Lemma \ref{lem_rank} that
\[
\dim (\Im (\partial)) \,  \leq \, \dim (\Im (d_{\ast, tH_{\rho}}))
\, =\,  \dim (\Im (d)),
\]
and hence $\dim (\Im (d) + \Im (\partial)) 
\leq 2 \dim (\Im (d)) = \dim (\Im (S)) $. 
Therefore $\Im (S) = \Im (d) + \Im (\partial)$. 
\qed

We are now ready to prove the disjointness of $d$ and $\partial$.

\bigskip
\noindent
{\bf Claim 3.} The two operators $d$ and $\partial$ are disjoint.

\bigskip
\noindent
{\bf Proof.} By Claim 2, both $\Im (d) \cap \Ker (\partial)$
and $\Im (\partial) \cap \Ker (d)$ are subspaces of $\Im (S)$.
They are clearly also subspaces of $\Ker (d) \cap \Ker (\partial)
\subset \Ker (S)$. Since $\Ker (S) \cap \Im (S) = 0$ by Claim 1,
we know that $\Im (d) \cap \Ker (\partial) = \Im (\partial) \cap 
\Ker (d) = 0$.
\qed

\bigskip
\begin{rem}
\label{rem_no-need-of-limit}
{\em
The fact that $\Im (\partial) \cap \Ker (d) = 0$
can be proved without using Theorem \ref{thm_limit-S}.
Indeed, suppose that $d \partial x = 0$ for some $x \in C$.
Then $\partial x \in \Ker (S)$. Thus $\phi (\partial x) = (1 + L_1 E)
\partial x \in \Ker (L)$. Since $\partial$ commutes with
$L_1E$ (see Remark \ref{rem_p-commute-E-L}),
we have $\partial (1 + L_1 E) x \in \Ker (L)$.
But we know from the usual Hodge theory
that $\Im (\partial) \cap \Ker (L) \subset \Im (L) \cap \Ker (L) = 0$.
Thus $(1 + L_1 E) x = 0$, or, $x = 0$. Therefore
$\Im (\partial) \cap \Ker (d) = 0$.
}
\end{rem}

\subsection{Poisson harmonic forms for $\pl$ as Hodge harmonic forms}
\label{sec_poi-hodge}

Consider now the bi-grading on $C = \Omega(K/T, \C)^K$ defined by
the complex structure $J$. Since $J$ is integrable, we
can write
\[
d \, = \, d^{'} \, + \, d^{''}
\]
where $d^{'}: C^{p,q} \lrw C^{p+1, q}$
and $d^{''}: C^{p,q} \lrw C^{p, q+1}$.
Since $C = \oplus_{p,q} C^{p,q}$ is an orthogonal
decomposition with respect to the Hermitian product
$\la \, , \, \ra_{\lambda}$, we have
\[
d_{\ast, \lambda} \, = \, d_{\ast, \lambda}^{'} \, + \, d_{\ast, 
\lambda}^{''},
\]
where $d^{'}_{\ast, \lambda}: C^{p,q} \lrw C^{p-1, q}$ and 
$d_{\ast, \lambda}^{''}: C^{p,q} \lrw C^{p, q-1}$ are
respectively the adjoint operators of $d^{'}$
and $d^{''}$ with respect to
$\la \, , \, \ra_{\lambda}$. Since $\ast_{\lambda}$
maps $C^{p,q}$ to $C^{{n \over 2} - q, {n \over 2} - p}$,
we know from $d_{\ast, \lambda} = (-1)^{p+q} 
\ast_{\lambda}^{-1} d \ast_{\lambda}$ on $C^{p,q}$
that 
\[
d_{\ast, \lambda}^{'} \, = \, (-1)^{p+q} 
\ast_{\lambda}^{-1} d^{''} \ast_{\lambda}, 
\hspace{.3in} {\rm and} \hspace{.3in}
d_{\ast, \lambda}^{''} \, = \, (-1)^{p+q} 
\ast_{\lambda}^{-1} d^{'} \ast_{\lambda}.
\]
Since $\pal = J^{-1} d_{\ast, \lambda} J$, we know that
$\pal (C^{p,q}) \subset C^{p-1,q} \oplus C^{p, q-1}$ (This fact 
also follows from Remark \ref{rem_simple-2}). Thus, 
we can write
\[
\pal \, = \, \partial_{\lambda}^{'} \, + \, 
\partial_{\lambda}^{''},
\]
where $\partial_{\lambda}^{'}: C^{p,q} \rightarrow
C^{p-1,q}$ and
$\partial_{\lambda}^{''}: C^{p,q} \rightarrow 
C^{p,q-1}$

The following proposition now follows immediately from the fact that
$J|_{C^{p,q}} = i^{p-q}\ {\rm id}$.

\begin{prop}
\label{prop_d-p}
The two decompositions $d_{\ast, \lambda} = d_{\ast, \lambda}^{'} 
+ d_{\ast, \lambda}^{''}$ and  $\pal = \partial_{\lambda}^{'} + 
\partial_{\lambda}^{''}$ are related by
\[
d_{\ast, \lambda}^{'} \, = \, - i \partial_{\lambda}^{'}
\hspace{.3in} {\rm and} \hspace{.3in}
d_{\ast, \lambda}^{''} \, = \, i \partial_{\lambda}^{''}.
\]
\end{prop}

The following is a corollary of Proposition \ref{prop_d-p}
and the fact that $d \pal + \pal d = 0$ on $C$:

\begin{cor}
\label{dor_degree-00}
On $C = \oplus_{p,q} C^{p,q}$, we have
\[
d^{'} d_{\ast, \lambda}^{''} \, + \, d_{\ast, \lambda}^{''} d^{'} \, = \, 0
\hspace{.3in} 
d^{''} d_{\ast, \lambda}^{'} \, + \, d_{\ast, \lambda}^{'} d^{''} \, = \, 0,
\]
and
\[
S_{\lambda} \, = \, 2(d^{'}  d_{\ast, \lambda}^{'} \, + \, 
 d_{\ast, \lambda}^{'} d^{'}) \, = \, 
2(d^{''}  d_{\ast, \lambda}^{''} \, + \,  
 d_{\ast, \lambda}^{''} d^{''})
\]
has bi-degree $(0,0)$. Consequently, the natural map
\[
\psi_{d, S_{\lambda}}: \, \, 
{\rm Ker} \,  (S_{\lambda}) \lrw H(C, d): \, \, 
\xi \Map [\xi]_d
\]
is an isomorphism of bi-degree $(0,0)$.
\end{cor}

We also immediately have

\begin{cor}
\label{cor_pal-harm-hodge-harm}
A $K$-invariant (complex valued) differential form $\xi$
on $K/T$ of pure bi-degree is harmonic with
respect to the Poisson structure $\pl$ (and a $K$-invariant
volume form $\mu_0$) if and only if
it is $(d, d_{\ast, \lambda})$-harmonic, i.e., it is Hodge harmonic
with respect to the Hermitian metric $h_{\lambda}$ on $K/T$.
\end{cor}

\begin{cor}
\label{cor_rep}
Every class in the de Rham cohomology of $K/T$ has
a unique representative that is Poisson harmonic
with respect to the Poisson structure $\pl$ (and
a $K$-invariant volume form $\mu_0$ on $K/T$).
\end{cor}

\begin{rem}
\label{rem_rep-ok-pl}
{\em 
The answer to Question \ref{ques_harm} is then ``yes" for the 
Poisson structure $\pi_{\lambda}$ and
 a $K$-invariant volume form $\mu_0$ on $K/T$.
}\end{rem}

\begin{rem}
\label{rem_compare}
{\em
It is interesting to compare the identities in this section
and the basic Kahler identities. See, for example, 
Corollary 4.10 in Chapter 5 of \cite{wells:complex}.
Compare also our operator $\partial_{\lambda} = J^{-1} d_{\ast, \lambda}J$
with the operator
$d_{c}^{*}$ on Page 191 in \cite{wells:complex}.
On the one hand, 
the Riemannian structure $g_{\lambda}$ and the complex structure $J$
we have here are not compatible to give a Kahler structure, but
on the other hand, the symplectic structure $\pl$ can be made 
into a Kahler form by a result of Ginzburg and Weinstein
\cite{g-w:poi}. It would be 
interesting to understand the connections here. We will do this 
elsewhere.
}
\end{rem}

\bigskip
To illustrate some special properties of the metrics $h_{\lambda}$
(see (\ref{eq_c-1}) and (\ref{eq_c-2}) below), we now
give another proof of Theorem
\ref{thm_limit-S}.

\bigskip
For notational simplicity, we will identify $C$ and $(\wedge \fn_{-}
\otimes \wedge \fn_{+})^T $ using the map $I$ and
regard the operators $d$ and $d_{\ast, \lambda}$
as operators on $(\wedge \fn_{-} \ot \wedge \fn_{+})^T$. We will
use $M_{\ast}$ and $M_{\ast, \lambda}$, for any operator $M$ on
$(\wedge \fn_{-} \ot \wedge \fn_{+})^T$, to denote the
adjoint operators of $M$ with respect to the Hermitian metrics
$h$ and $h_{\lambda}$ respectively. For each $\alpha > 0$, set
\[
c_{\alpha} (\lambda) \, = \, {e^{\alpha(\lambda)} \over
e^{2\alpha (\lambda)} -1}.
\]
Let $A_{\lambda} \in \End (C)$ be the multi-linear extension of
the operator on $\fn_{-} \oplus \fn_{+}$ defined by
\[
A_{\lambda} \ea \, = \, c_{\alpha}(\lambda) \ea, \hspace{.2in}
A_{\lambda} \eb \, = \, c_{\alpha}(\lambda) \eb.
\]
Then $h_{\lambda}(\xi, \eta) = h(A_{\lambda} \xi, \eta)$ for all
$\xi, \eta \in (\wedge \fn_{-} \ot \wedge \fn_{+})^T$, and
\[
M_{\ast, \lambda} \, = \, A_{\lambda}^{-1} M_{\ast} A_{\lambda}
\]
for any $M \in \End(C)$. In particular, we have
\[
d_{\ast, \lambda}^{'} \, = \,  A_{\lambda}^{-1} d_{\ast}^{'} A_{\lambda}.
\]
But, by identifying $\fn_{-}$ with $\fn_{+}^{*}$
via the Killing form of $\fg$, we can identify  $d^{'}$ with
the Chevalley-Eilenberg coboundary operator for the Lie
algebra $\fn_{+}$ with coefficients in $\wedge \fn_{+}$ as the
adjoint module (see, for example, (3.2.1) in \cite{ko:63}). Thus,
we can further decompose $d^{'}$ as
\[
d^{'} \, = \, d_{{\frak n}_{+}} \ot 1 \, + \, \kappa,
\]
where $d_{{\frak n}_{+}}$ is the Chevalley-Eilenberg coboundary operator
for the Lie algebra $\fn_{+}$ with trivial
coefficients, and
\[
\kappa: \, \, (\wedge^p \fn_{-} \ot \wedge^q \fn_{+})^T
\lrw (\wedge^{p+1} \fn_{-} \ot \wedge^q \fn_{+})^T
\]
is given by
\[
\kappa (\xi \ot X) \, = \, (-1)^p \xi \wedge \sum_{\alpha > 0}
\eb \, \ot \, {\rm ad}_{E_{\alpha}} X.
\]
Thus we have
\[
d^{'}_{\ast, \lambda} \, = \, (d_{{\frak n}_{+}} \ot 1)_{\ast, \lambda}
\, + \, \kappa_{\ast, \lambda} \, = \,
A_{\lambda}^{-1} (d_{{\frak n}_{+}} \ot 1)_{\ast} A_{\lambda}
\, + \, A_{\lambda}^{-1} \kappa_{\ast} A_{\lambda}.
\]
It is easy to see that $(d_{{\frak n}_{+}} \ot 1)_{\ast} =
- b_{{\frak n}_{-}}
\ot 1$, where $b_{{\frak n}_{-}}$ is the Chevalley-Eilenberg
boundary operator for the Lie algebra $\fn_{-}$, and that
$\kappa_{\ast}$ is given by
\[
\kappa_{\ast} (\eta \ot Y) \, = \, \sum_{\alpha >0} i_{E_{\alpha}} \eta
\ot {\rm ad}_{-E_{\alpha}}^{*} Y.
\]
Now an easy calculation using the fact
\begin{equation}
\label{eq_c-1}
\lim_{t \rightarrow +\infty} {c_{\alpha} c_{\beta} \over c_{\alpha +
\beta}} (\lambda + tH_{\rho}) \, = \, 1
\end{equation}
if $\alpha + \beta $ is  a   root
shows that
\[
\lim_{t \rightarrow +\infty} A^{-1}_{\lambda + tH_{\rho}}
(b_{{\frak n}_{-}} \ot 1) A_{\lambda + tH_{\rho}} \, = \, 
b_{{\frak n}_{-}} \ot
1.
\]
Similarly, using the fact that
\begin{equation}
\label{eq_c-2}
\lim_{t \rightarrow +\infty} {c_{\alpha} c_{\beta} \over 
c_{\beta - \alpha
}} (\lambda + tH_{\rho}) \, = \, 0
\end{equation}
if $\beta - \alpha$ is a positive root,
we can show that
\[
\lim_{t \rightarrow +\infty} A^{-1}_{\lambda + tH_{\rho}}
\kappa_{\ast} A_{\lambda + tH_{\rho}} \, = 0.
\]
Thus
\[
\lim_{t \rightarrow +\infty} d^{'}_{\ast, \lambda + tH_{\rho}} \, = \,
- (b_{{\frak n}_{-}} \ot 1) \, = \, \partial^{'}.
\]
Similarly, we can show that
\[
\lim_{t \rightarrow +\infty} d^{''}_{\ast, \lambda+ tH_{\rho}} \, = \,
\partial^{''}.
\]
Hence $ \lim_{t \rightarrow +\infty} d_{\ast, \lambda + 
tH_{\rho}} = \partial$ and
$ \lim_{t \rightarrow +\infty} S_{\lambda +tH_{\rho}} = S$.

\subsection{The Kostant harmonic forms as limits of Hodge harmonic forms}
\label{sec_limit-indeed}
 
As we have seen from Example \ref{exam_not-adjoint}, the
Kostant operator $\partial$ is not the adjoint of the de Rham $d$
with respect to any Hermitian metric on $K/T$, so
the Kostant harmonic forms $s^w, w \in W$, are not harmonic in the
sense of Hodge. However, we will show by using Theorem
\ref{thm_limit-S} that they are limits of Hodge harmonic
forms with respect to the family of Hermitian metrics
$h_{\lambda}$.
 
\bigskip
Recall from Section \ref{sec_another-proof} that
$C  =  \Ker (S)  +  \Im (S)$
is a direct sum decomposition. Consider the linear map
\[
\psi_{S, S_{\lambda}}: \, \,  \Ker (S_{\lambda}) \lrw \Ker (S):\, \,
x \Map x_1, \hspace{.2in} x \in \Ker (S_{\lambda})
\]
if $x = x_1 + x_2$ with $x_1 \in \Ker (S)$ and $x_2 \in \Im (S)$.
Since
\[
\Ker (\psi_{S, S_{\lambda}}) \, \subset \,
\Ker (S_{\lambda}) \cap \Im (S)
\, =  \, \Ker (S_{\lambda}) \cap \Ker (d) \cap \Im (S) \, \subset \,
\Ker (S_{\lambda}) \cap \Im (d) \, = \, 0,
\]
the map $\psi_{S, S_{\lambda}}$ is injective. But since both
$\Ker (S)$ and $\Ker (S_{\lambda})$ have the same dimension as
$H^{\bullet}(K/T, \C)$, the map $\psi_{S, S_{\lambda}}$
is an isomorphism.
 
\bigskip
\begin{dfn}
\label{dfn_swlambda}
{\em
For $w \in W$, we define the Hodge harmonic form $s^{w}_{\lambda}$
(with respect to the Hermitian metric $h_{\lambda}$) to be
\[
s^{w}_{\lambda} \, = \, (\psi_{S, S_{\lambda}})^{-1} s^w \,
\in \, \Ker (S_{\lambda}).
\]
}
\end{dfn}
 
\bigskip
\begin{thm}
\label{thm_main-4}
For each $w \in W$, we have
\[
\lim_{t \rightarrow + \infty} s^{w}_{\lambda + t H_{\rho}} \, = \, s^w.
\]
\end{thm}
 
\bigskip
\noindent
{\bf Proof.} Since $C = \Ker (S) + \Im (S)$ is a direct sum,
we have the Green's operator $S_1$ of $S$, namely $S_1 \in \End(C)$
is such that $S_1|_{{\rm Ker (S)}} = 0$ and
$S_1|_{{\rm Im}(S)} = (S|_{{\rm Im}(S)})^{-1}$. Set
$F_{\lambda} = S_{\lambda} - S$. Then, for each $w \in W$,
\[
s^w \,  = \,  s^{w}_{\lambda} \, - \, S_1 S (s^{w}_{\lambda}) 
\,  =\,  s^{w}_{\lambda} \, + \, S_1 F_{\lambda} (s^{w}_{\lambda}) 
\,  = \,  (1 \, + \, S_1 F_{\lambda}) (s^{w}_{\lambda}).
\]
Since $F_{\lambda + t H_{\rho}}$ goes to $0$ as $t \rightarrow +\infty$,
the operator $1 \, + \, S_1 F_{\lambda + tH_{\rho}}$
 is invertible for $t$ large
enough, and thus
\[
\lim_{t \rightarrow +\infty} s^{w}_{\lambda + tH_{\rho}}
\, = \, \lim_{t \rightarrow +\infty} (1 \, + \, S_1 
F_{\lambda + tH_{\rho}})^{-1}
s^w \, = \, s^w.
\]
\qed

\section{Appendix: The Schouten bracket}
\label{sec_appendix}

Let $P$ be a manifold of dimension $n$.
The Schouten bracket $[ \, , \, ]$  is the graded Lie bracket
on the space  $\chi^{\bullet}(P)$
of multi-vector fields on $P$ that
is characterized as follows \cite{kz:schouten}:

1) for any smooth function $f \in \chi^0(P) =
C^{\infty}(P)$ and any ($1$-)vector field
$X \in \chi^1(P)$,
\[
[X, ~ f] ~ = ~ X(f) \in C^{\infty}(P);
\]

2) for any ($1$-)vector fields $X$ and $Y$, the Schouten bracket
$[X, ~ Y]$ is the usual Lie bracket between $X$ and $Y$.
 
3) the bracket between two general multi-vector fields
is obtained according to the following two rules:
\beqa
& & [X, Y] =  - (-1)^{(|\tx|-1)(| \ty|-1)} [Y, X];\\
& & [X, Y \wedge Z] = [X, Y] \wedge Z + (-1)^{(|\tx|-1)|\ty|} Y \wedge
 [X, Z].
\eeqa
It satisfies the graded Jacobi identity:
\[
(-1)^{(|X|-1)(|Z|-1)}[X, \, [Y, \, Z]] \, + \, c.p.(X, Y, Z) \, = \, 0,
\]
where $c.p.(X, Y, Z)$ means the cyclic permutation between
$X, Y$ and $Z$.
Explicitly, if $\alpha$ is a form on $P$ of
degree $|X| + |Y| -1$,
we have,
\begin{equation}
\label{eq_sch-explicit-on-p}
(\alpha, ~ [X, Y]) ~ = ~ (-1)^{(|\tx|-1)(|\ty|-1)} i_{\tx}
di_{\ty} \alpha ~ - ~ i_{\ty} d i_{\tx} \alpha ~ + ~
(-1)^{|\tx|} i_{\tx \wedge \ty} d \alpha.
\end{equation}
For $X \in \chi^{|\tx|}(P)$, introduce the operator
\[
\partial_{\tx} \, = \, i_{\tx} d \, - \, (-1)^{|\tx|}
d i_{\tx}: \, \,
\Omega^k(P) \lrw \Omega^{k-|\tx| + 1}(P).
\]
Then it follows from the definition of the Schouten
bracket that
\begin{eqnarray}
\label{eq_i-wedge}
& & i_{\tx \wedge \ty} ~ = ~ i_{\ty} i_{\tx},\\
\label{eq_i-bra}
& & i_{[X, Y]} ~ = ~ (-1)^{(|X|-1)(|Y|-1)}(
\partial_X i_Y ~ - ~ (-1)^{(|X|-1)|Y|} i_Y \partial_X),\\
\label{eq_p-bra}
& & \partial_{[X, Y]} ~ = ~(-1)^{(|X|-1)(|Y|-1)}( \partial_X \partial_Y ~ - ~
(-1)^{(|X|-1)(|Y|-1)} \partial_Y \partial_X).
\end{eqnarray}
We can think of the operators $i_X$ and $\partial_X$ as defining a
``right" representation of
the Gerstenhaber algebra
$\chi(P)$ \cite{ks:gerst} \cite{xu:modular} on the graded vector space 
$\Omega(P)$.

\end{document}